\documentclass[12pt]{article}

\usepackage{amsfonts}
\usepackage{amsbsy}
\usepackage{graphicx}

\def\blacksquare{\vrule height .9ex width .8ex depth -.1ex}

\newcommand{\ih}{\'{\i}}
\newcommand{\eh}{\hspace{.06in}}

\newcommand{\C}{\mathbb{C}}
\newcommand{\R}{\mathbb{R}}
\newcommand{\Z}{\mathbb{Z}}
\newcommand{\cH}{\mathcal{H}}
\newcommand{\cR}{\mathcal{R}}
\newcommand{\hC}{\hat{\C}}
\newcommand{\af}{\alpha}
\newcommand{\bt}{\beta}
\newcommand{\g}{\gamma}
\newcommand{\ld}{\lambda}
\newcommand{\G}{\Gamma}
\newcommand{\m}{{\frac{1}{2}}}
\newcommand{\pd}{\pi/2}
\newcommand{\tw}{\tilde{w}}
\newcommand{\tT}{\tilde{T}}
\newcommand{\impl}{\Rightarrow}
\newcommand{\see}{\eh{\Leftrightarrow}\eh}
\newcommand{\ds}{\displaystyle}
\newcommand{\deh}{\partial}
\newcommand{\ovl}{\overline}
\newcommand{\sen}{\rm\,sen\,}
\newcommand{\BE}{\begin{equation}}
\newcommand{\EE}{\end{equation}}
\newcommand{\fns}{\large}

\begin{document}
\thispagestyle{empty}
\centerline{\bf\huge FUN\c C\~OES \ EL\'ITICAS }
\bigskip
\bigskip
\bigskip

\centerline{\bf\Large K{\fns ELLY} R{\fns OBERTA} M{\fns AZZUTTI} L{\fns \"UBECK}}

\centerline{\large Universidade Estadual do Oeste do Paran\'a - UNIOESTE, Brasil}\bigskip
\bigskip

\centerline{\bf\Large V{\fns AL\'ERIO} R{\fns AMOS} B{\fns ATISTA}}

\centerline{\large Universidade Federal do ABC - UFABC, Brasil}
\eject
\thispagestyle{empty}
\bigskip
\bigskip
\ \\
{\bf Resumo:} Em 1989, H.Karcher re-escreveu a teoria de fun\c c\~oes el\ih ticas utilizando uma abordagem bem mais geom\'etrica do que anal\ih tica. Com isso, obteve um \'otimo controle sobre o comportamento e valores-imagem destas fun\c c\~oes, permitindo sua larga aplica\c c\~ao em superf\ih cies m\ih nimas. Este trabalho destina-se a apresentar a teoria de fun\c c\~oes el\ih ticas segundo Karcher, bem como algumas aplica\c c\~oes desta em superf\ih cies m\ih nimas.
\eject
\thispagestyle{empty}
\ \\
{\huge\bf \'Indice}
\bigskip
\bigskip
\ \\
{\bf Cap\ih tulo 1}
\\
{\bf Fun\c c\~oes El\ih ticas}\hfill 1
\\
1. Introdu\c c\~ao\dotfill 1
\\
2. A rela\c c\~ao entre $P$ e $\wp$\dotfill 3
\\
3. Uma observa\c c\~ao importante\dotfill 5
\\
\\
{\bf Cap\ih tulo 2}
\\
{\bf A fun\c c\~ao $\wp$-Weierstra\ss \ Sim\'etrica}\hfill 6
\\
1. Resultados b\'asicos\dotfill 6
\\
2. A fun\c c\~ao $\wp$-Weierstra\ss \ sim\'etrica\dotfill 9
\\
3. A equa\c c\~ao alg\'ebrica do toro\dotfill 14
\\
4. A constante $c$\dotfill 16
\\
\\
{\bf Cap\ih tulo 3}
\\
{\bf A fun\c c\~ao $\g$}\hfill 19
\\
\\
{\bf Cap\ih tulo 4}
\\
{\bf Aplica\c c\~oes}\hfill 25
\\
1. Resultados cl\'assicos\dotfill 25
\\
2. Exemplo de superf\ih cie m\ih nima\dotfill 27
\\
\\
{\bf Refer\^encias}\hfill 39
\eject
\count0=1
\ \\
{\huge\bf Cap\ih tulo 1}
\bigskip
\bigskip
\ \\
{\huge\bf Fun\c c\~oes El\ih ticas}
\bigskip
\bigskip
\bigskip
\ \\
{\bf 1. Introdu\c c\~ao}
\\

Dizemos que $F:\C\to\hC$ \'e uma {\it fun\c c\~ao el\ih tica} se ela for {\it duplamente peri\'odica} \ e {\it meromorfa}. Ou seja, em primeiro lugar, existe uma {\it base} $\{w_1,w_2\}$ de $\C=\R^2$, como espa\c co vetorial real, de modo que $F(z+w_1)=F(z+w_2)=F(z)$, $\forall\,z\in\C$. Neste sentido, subentende-se que $w_1$ e $w_2$ s\~ao os {\it per\ih odos} de $F$, donde $\forall\,w$ tal que $F(\cdot+w)\equiv F$ implica $w=nw_1+mw_2$ para certos $m,n\in\Z$. 
\\

Note que isso nos impede de estabelecer $\{w_1,w_2\}$ no caso de fun\c c\~oes cons- tantes. Assim, em segundo lugar, os conjuntos $A=F^{-1}(\C)$ e $B=F^{-1}(\hC^*)$ s\~ao ambos n\~ao-vazios, e desta forma exigimos que $F|_A$ e $1/F|_B$ sejam fun\c c\~oes {\it complexas anal\ih ticas}. Essa propriedade traduz o que quizemos dizer por ``$F$ meromorfa''. 
\\

Quanto \`as fun\c c\~oes constantes, iremos tamb\'em consider\'a-las el\ih ticas por extens\~ao. Isso porque, deste modo, podemos dizer que {\it toda} fun\c c\~ao el\ih tica est\'a definida num toro.
\\

De fato, dada uma base $\{w_1,w_2\}\in\C$, esta determina um {\it reticulado} $\G=\{nw_1+mw_2:n,m\in\Z\}$ e a rela\c c\~ao $z\sim\zeta\see z-\zeta\in w_1\Z+w_2\Z$, que \'e de {\it equival\^encia}. As propriedades reflexiva e sim\'etrica s\~ao \'obvias, enquanto que $z\sim\zeta\sim w\impl\zeta-w=Nw_1+Mw_2$, $z-\zeta=nw_1+mw_2$, donde $z-w=(n+N)w_1+(m+M)w_2\in\G$. Com isso temos as {\it classes de equival\^encia} $[z]$, onde $z$ \'e sempre um representante qualquer, e o conjunto de tais classes \'e o {\it toro} $T=\C/\G$. Podemos ent\~ao introduzir a {\it fun\c c\~ao quociente} $p:\C\to\C/\G$, dada por $p(z)=[z]$.
\\

Note uma sutileza: apesar de $\G=w_1\Z+w_2\Z$, usamos sempre $w_1\Z+w_2\Z$ na defini\c c\~ao das classes de equival\^encia. Com isso, os reticulados $w+\G$, $\forall\,w\in\Z$, s\~ao todos distintos, mas definem exatamente o {\it mesmo} toro.\\
\\
{\bf Defini\c c\~ao 1.1.} Sejam $X$ e $Y$ espa\c cos topol\'ogicos. Uma aplica\c c\~ao $r:Y\to X$ \'e dita {\it recobrimento} se qualquer $x\in X$ admite uma vizinhan\c ca aberta $U$ tal que $r^{-1}(U)=\dot\cup_\ld V_\ld$, onde os $V_\ld$ s\~ao abertos de $Y$, dois a dois disjuntos, e $r|_{V_\ld}$ \'e um {\it homeomorfismo} $\forall\,\ld$.
\\

De acordo com a Defini\c c\~ao 1.1, vemos facilmente que $p:\C\to T$ \'e um recobrimento. Assim, qualquer fun\c c\~ao el\ih tica $F:\C\to\hC$ com per\ih odos $w_{1,2}$ pode ser vista como $f:T\to\hC$ simplesmente por $f([z])=F(z)$, e \'e claro que isso inclui as fun\c c\~oes constantes. 
\\

Mais tarde veremos que toda fun\c c\~ao el\ih tica \'e a inversa de uma {\it integral el\ih tica}, que define-se como do tipo $\int_c^x R(t,\sqrt{q(t)})dt$, onde $R$ \'e fun\c c\~ao racional, $q$ \'e um polin\^omio de grau 3 ou 4 sem ra\ih zes m\'ultiplas, e $c$ \'e uma constante. O exemplo que originou esta denomina\c c\~ao foi $\ds\int_0^x\frac{\sqrt{(1-t^2)(1+\kappa^2t^2)}}{1-t^2}dt$, repre- sentando o comprimento de arco da {\it elipse} com raios principais $1$ e $\kappa$.  
\\

Uma das primeiras fun\c c\~oes el\ih ticas, cujo estudo encontramos com detalhes em v\'arias literaturas de an\'alise complexa, \'e a fun\c c\~ao ``P de Weierstra\ss'', simbolizada por $\wp$. Entretanto, aqui reservaremos este s\ih mbolo para a ``P de Weierstra\ss \ sim\'etrica'', a ser introduziada no Cap\ih tulo 2. Quanto \`a ``P de Weierstra\ss \ cl\'assica'', esta denotaremos simplesmente por $P$. Sua defini\c c\~ao original \'e dada por uma soma sobre $\G^*=\G\setminus\{0\}$, qual seja 
\BE
   P(z)=\frac{1}{z^2}+\sum_{w\in\G^*}\biggl[\frac{1}{(z-w)^2}-\frac{1}{w^2}\biggl].
\EE

Prova-se que (1) est\'a bem definida, isto \'e, independe de como se enumera os elementos de $\G$. Al\'em disso, temos que $P$ satisfaz \`a Equa\c c\~ao Diferencial 
\[
   P^{\prime 2}=4P^3-\af P-\bt,
\]
Onde $\af,\bt$ s\~ao complexos que dependem de $w_1/w_2$. Dada 
\[
   y(x)=\int_x^\infty\frac{dt}{\sqrt{4t^3-\af t-\bt}},
\]
mostra-se ainda que $x=P(y)$. Estes fatos deixaremos aqui apenas como coment\'ario, pois percebe-se que a defini\c c\~ao de $P$ \'e fortemente anal\ih tica, e de dif\ih cil manipula\c c\~ao. No Cap\ih tulo 2 iremos introduzir a func\~ao $\wp$, cuja invers\~ao $1/\wp$ \'e apenas uma transforma\c c\~ao linear afim de $P$, por\'em muito mais simples de ser estudada.  
\\
\\
{\bf 2. A rela\c c\~ao entre $P$ e $\wp$}
\\

Para deduzir a rela\c c\~ao entre essas duas fun\c c\~oes el\ih ticas, \'e mais f\'acil considerar a base $\{2w_1,2w_2\}$ e o reticulado $G=\{(2n+1)w_1+(2m+1)w_2:n,m\in\Z\}$. Tomando $W_{1,2}=2w_{1,2}$ e $\G=\{nW_1+mW_2:n,m\in\Z\}$, as respectivas bases geram as mesmas classes de equival\^encia para $\C/\G$ e $\C/G$, donde a igualdade entre esses dois toros.  
\\

A Figura 1.1 representa um paralelogramo fundamental do reticulado $G$. J\'a vimos que $P$ tem um {\it p\'olo} de ordem dois na origem, e (1) implica que este \'e seu \'unico p\'olo no toro. Assim, $P(w_1)=e_1$, $P(w_2)=e_2$ e $P(w_1+w_2)=e_3$, todos valores em $\C$.   
\\
\begin{figure}[!ht]
\centering
\includegraphics[scale=0.9]{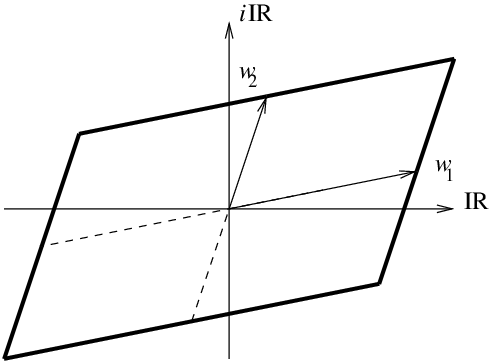}
\centerline{Figura 1.1: O toro $T$.}
\end{figure}

Sem entrar em detalhes que ser\~ao vistos no Cap\ih tulo 2, vamos aqui somente adiantar que $\wp(0)=0^2$ e $\wp(w_1+w_2)=\infty^2$, sendo estes seus \'unicos zero e p\'olo em $T$. Como ambas s\~ao meromorfas, temos as seguintes {\it s\'eries de Laurent} numa vizinha\c ca da origem:
\[
   P(z)=\frac{1}{z^2}+\frac{a_{-1}}{z}+a_0+a_1z+a_2z^2+\cdots\eh\eh{\rm e}
\]    
\[
   \frac{1}{\wp(z)}=\frac{b_{-2}}{z^2}+\frac{b_{-1}}{z}+b_0+b_1z+b_2z^2+\cdots
\]   

Assim, a func\~ao $f=b_{-2}P-1/\wp$ \'e meromorfa no toro, mas talvez n\~ao possua nenhum p\'olo, o que chamaremos de {\it holomorfa}, termo que ser\'a melhor explicado na Se\c c\~ao 3. Neste caso, se $f$ n\~ao for constante, a teoria de an\'alise complexa garante que ela \'e aberta. Mas sendo $T$ compacto, e $f$ cont\ih nua, ent\~ao $f(T)\ne\emptyset$ \'e ao mesmo tempo um aberto e um fechado de $\hC$, donde $f(T)=\hC$. Ou seja, $f$ deveria assumir p\'olo. A \'unica possibilidade \'e que seja um p\'olo simples na origem, donde
\BE
   f(z)=\frac{c_{-1}}{z}+c_0+c_1z+c_2z^2+\cdots
\EE   
e $c_{-1}=b_{-2}a_{-1}-b_{-1}\ne 0$. Agora usaremos o seguinte resultado da Teoria de Superf\ih cies de Riemman:
\\
\\
{\bf Teorema 1.1} \it Seja $f:T\to\hC$ meromorfa n\~ao constante. Ent\~ao $f^{-1}(z)$ e $f^{-1}(w)$ t\^em sempre o mesmo n\'umero de elementos, $\forall\,z,w\in\hC$, cada qual contando com sua ``multiplicidade''. Ou seja, tomando a correspondente $F:\C\to\hC$ e $\zeta\in F^{-1}(z)$, $z\in\C$, \'e ``a primeira ordem $m\ge 1$ de deriva\c c\~ao $F^{(m)}$ que n\~ao se anula em $\zeta$'' (para $z=\infty$ considera-se $(1/F)^{(m)}$).\rm  
\\

Chamamos {\it grau} de $f$ a cardinalidade de $f^{-1}(z)$, $\forall\,z\in\hC$, que denotamos por deg$(f)$. Pelo Teorema 1.1, vemos que a $f$ encontrada em (2) \'e uma bije\c c\~ao de $T$ em $\hC$, em particular cont\ih nua. Como $f$ \'e aberta, sua inversa \'e tamb\'em cont\ih nua. Ou seja, $f$ \'e um {\it homeomorfismo} entre $T$ e $\hC$. Ocorre que $\hC$ \'e {\it simplesmente conexo}, mas isso n\~ao vale para $T$. Assim, chegamos a uma contradi\c c\~ao, que deveu-se ao fato de supormos $f$ n\~ao-constante.
\\

Ent\~ao, como $f$ \'e uma constante, temos $P=a/\wp+b$ com $a,b\in\C$ e $a\ne 0$.
\\
\\
{\bf 3. Uma observa\c c\~ao importante}
\\
 
Na introdu\c c\~ao deste cap\ih tulo, discutimos o que se entende por uma fun\c c\~ao $F$ ``meromorfa''. Mas este conceito tamb\'em compreende o fato de que Dom$(F)$ \'e um aberto de $\C$ e $F$ pode assumir o valor $\infty$. Quando o dom\ih nio \'e uma superf\ih cie compacta, ou se Im$(F)\subset\C$, ent\~ao dizemos que \'e uma fun\c c\~ao {\it holomorfa}. Al\'em disso, se existe a inversa $F^{-1}$, dizemos que \'e {\it biholomorfa}. Todos esses termos admitem o prefixo {\it anti-} quando aplicados \`a conjugada $\ovl{F}$. Por exemplo, $F$ anti-holomorfa equivale a $\ovl{F}$ holomorfa, etc.
\eject
\ \\
{\huge\bf Cap\ih tulo 2}
\bigskip
\bigskip
\ \\
{\huge\bf A Fun\c c\~ao $\wp$-Weierstra\ss \ Sim\'etrica}
\bigskip
\bigskip
\bigskip

Em compara\c c\~ao com a $P$-Weierstra\ss \ cl\'assica, a $\wp$-Weierstra\ss \ sim\'etrica tem as seguintes vantagens: \'e muito f\'acil de ser manipulada geometricamente, e tamb\'em oferece bem mais informa\c c\~oes a respeito de seus valores-imagem ao longo do toro. Como acabamos de ver na se\c c\~ao anterior, ambas guardam uma rela\c c\~ao alg\'ebrica muito simples: $P=a/\wp+b$ com $a,b\in\C$ e $a\ne 0$. 
\\
\\
{\bf 1. Resultados b\'asicos}
\\

Neste cap\ih tulo, trabalharemos com o cenceito de {\it superf\ih cie de Riemann}, que \'e um par $(X,\Sigma)$, onde $X$ \'e uma variedade 2-dimensional conexa
e $\Sigma$ \'e uma {\it estrutura complexa} em $X$. Al\'em dos abertos conexos de $\hC$, no Cap\ih ulo 1 vimos o toro, que tamb\'em \'e exemplo de superf\ih cie de Riemann. De fato, como $p:\C\to\C/\G$ \'e recobrimento, podemos induzir a estrutura complexa de $\C$ em $T$ do seguinte modo: seja $V\subset\C$ um subconjunto aberto tal que quaisquer dois de seus pontos n\~ao s\~ao equivalentes m\'odulo $\G$. Ent\~ao $U=p(V)$ \'e aberto e $p|_V:V\to U$ \'e um homeomorfismo. Sua {\it inversa local} $(p|_V)^{-1}:U\to V$ \'e uma {\it carta}, e a cole\c c\~ao de todas as cartas assim formadas constitui a estrutura $\Sigma$ de que falamos anteriormente. 
\\

Vamos introduzir agora resultados importantes que n\~ao est\~ao estritamente relacionados com fun\c c\~oes el\ih ticas. Estes resultados constituem os argumentos fundamentais para podermos analisar o comportamento das fun\c c\~oes el\ih ticas. Primeiramente, precisamos da seguinte defini\c c\~ao:
\\
\\
{\bf Defini\c c\~ao 2.1.} Uma {\it involu\c c\~ao} \'e uma aplica\c c\~ao cont\ih nua $I:S\to S$ que satisfaz $I\circ I=id$ (a identidade em $S$). Quando $S$ \'e uma superf\ih cie compacta, a involu\c c\~ao \'e chamada {\it hiperel\ih tica} se $S/I$ \'e homeomorfa a $S^2$.
\\

Da defini\c c\~ao acima, \'e imediato ver que toda involu\c c\~ao \'e uma bije\c c\~ao. No pr\'oximo teorema, mencionamos as {\it transforma\c c\~oes de M\"obius}, que s\~ao aplica\c c\~oes $T:\hC\to\hC$ dadas por $T(z)=\ds\frac{az+b}{cz+d}$ com $a,b,c,d$ constantes em $\C$ e $ad\ne bc$. 
\\

Um resultado conhecido da an\'alise complexa \'e que {\it todo} biholomorfismo de $\hC$ em $\hC$ \'e caracterizado por uma transforma\c c\~ao de M\"obius. Al\'em disso, uma tal transforma\c c\~ao sempre leva circunfer\^encias em circunfer\^encias (de $\hC$). Vejamos agora o que ocorre com as involu\c c\~oes:
\\
\\
{\bf Teorema 2.1.} \it Toda involu\c c\~ao holomorfa ou anti-holomorfa em $\hC$ \'e dada por uma transforma\c c\~ao de M\"obius $M$ ou sua conjugada $\ovl{M}$.\rm\\
\\
\underline{Prova}: Tome uma involu\c c\~ao $I:\hC\to\hC$. Se \'e holomorfa, ent\~ao \'e biholomorfa e assim uma transforma\c c\~ao de M\"obius. Se \'e anti-holomorfa, ent\~ao $\ovl{I}$ \'e biholomorfa, logo $I$ \'e a conjugada de uma transforma\c c\~ao de M\"obius\hfill c.q.d.
\\

Note que a rec\ih proca do teorema acima {\it n\~ao} \'e v\'alida, pois $z\to 2z$ n\~ao \'e uma involu\c c\~ao.
\\
\\
{\bf Teorema 2.2.} \it Sejam $S$ e $R$ superf\ih cies de Riemann, $I:S\to S$ uma involu\c c\~ao e $f:S\to R$ uma fun\c c\~ao cont\ih nua, aberta e sobrejetora. Ent\~ao, existe uma \'unica involu\c c\~ao $J:R\to R$ tal que $J\circ f=f\circ I$ se, e somente se, sempre que $f(x)=f(y)$ temos $f\circ I(x)=f\circ I(y)$.\rm
\\
\\
\underline{Prova}: As hip\'oteses do teorema foram formuladas para garantir a comutatividade do diagrama:
\begin{figure}[h]
\center
\includegraphics[scale=1]{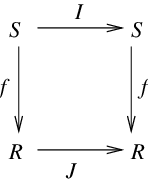}\\
\end{figure}

Admitimos que $J\circ f=f\circ I$. Se $f(x)=f(y)$ ent\~ao 
\[
   f\circ I(x)=J\circ f(x)=J\circ f(y)=f\circ I(y).
\]

Suponhamos agora que $f\circ I(x)=f\circ I(y)$ quando $f(x)=f(y)$. Para cada $z\in R$ definimos $J(z):=f\circ I(x)$, para algum $x\in S$ tal que $f(x)=z$ (este elemento existe pois $f$ \'e sobrejetora). Por hip\'otese $J:S\to S$ est\'a bem definida. Al\'em disso, $J\circ f=f\circ I.$ A fun\c c\~ao $J$ \'e cont\ih nua por causa do seguinte argumento: para qualquer subconjunto aberto $U\subset R$ temos que $I^{-1}(f^{-1}(U))$ \'e aberto em $S$. Mas $I^{-1}(f^{-1}(U))=f^{-1}(J^{-1}(U))$, e conseq\"uentemente $J^{-1}(U)$ \'e aberto em $R$. Temos ainda
\[
   J\circ J\circ f=J\circ f\circ I=f\circ I\circ I=f.
\]

Assim, $J$ \'e uma involu\c c\~ao\hfill c.q.d.
\\
\\
{\bf Observa\c c\~ao.} Uma vers\~ao mais geral do Teorema 2.2 \'e a seguinte: \it Sejam $S$ e $R$ superf\ih cies de Riemann, $I:S\to S$ um biholomorfismo e $p,q:S\to R$ recobrimentos. Ent\~ao, $(\exists\,!\,\,J:R\to R$ biholomorfa tal que $J\circ p=q\circ I)\see$ $(q\circ I(x)=q\circ I(y)\impl p(x)=p(y))$.\rm 
\\

\'E imediato concluir que dado $F\subset R$ temos $J(F)=F$ se, e somente se, $I(f^{-1}(F))=f^{-1}(F)$.
\\
\\
\\
{\bf 2. A fun\c c\~ao $\wp$-Weierstra\ss \  sim\'etrica}
\\

Como vimos no Cap\ih tulo 1, qualquer toro $T$ \'e o quociente de $\C$ por algum reticulado $G\subset\C$. Podemos descrever isto por meio de n\'umeros complexos n\~ao nulos $w_1$ e $w_2,$ com $\ds\frac{w_1}{w_2}\notin\R$. Assim, $G=\{(2n+1)w_1+(2m+1)w_2: n,m\in\Z\}$ e $T=\C/G$.
\\
\begin{figure}[h]
\center
\includegraphics[scale=0.75]{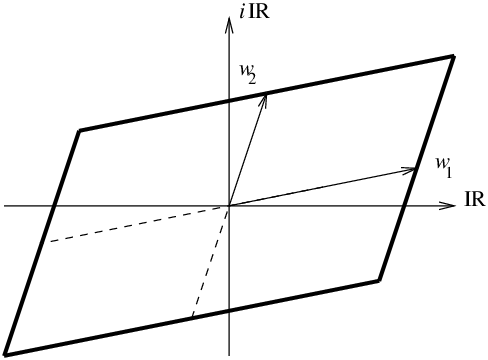}
\centerline{Figura 2.1: O Toro $T$.}
\end{figure}

No toro $T$ a involu\c c\~ao $I(z)=-z$ tem exatamente quatro pontos fixos: $0,w_1,w_2$ e $w_1+w_2$. Usando a {\it f\'ormula de Euler-Poincar\'e} obtemos:
\[
   \chi(T/I)=\frac{\chi(T)}{2}+\frac{4}{2}=2.
\]

Logo, o g\^enero de $T/I$ \'e zero e pelo {\it teorema de Koebe} existe um biholomorfismo $\mathcal{B}:T/I\to\hC$. A menos de uma transforma\c c\~ao de M\"obius a fun\c c\~ao $\mathcal{B}$ est\'a bem definida.
\\
\\
{\bf Defini\c c\~ao 2.2.} A {\it fun\c c\~ao $\wp$-Weierstra\ss \  sim\'etrica} \'e a composta $\mathcal{B}\circ(\cdot/I):T\to\hC$ tal que
\[
   \wp(0)=0,\eh\eh\wp(w_1+w_2)=\infty\eh\eh{\rm e}\eh\eh
   \wp\biggl(\frac{w_1+w_2}{2}\biggl)=i.
\]

Devido \`a fun\c c\~ao quociente $\cdot/I:T\to T/I$, \'e imediato concluir que $\wp(-z)=\wp(z),$ para qualquer $z\in T.$ Logo $\wp\biggl(-\frac{\ds w_1+w_2}{\ds 2}\biggl)=i.$ Al\'em disso, deg($\wp$)=2 e seus pontos de ramo s\~ao exatamente os quatro pontos fixos de $I$.
\\

Da Defini\c c\~ao 2.1 vemos que $I$ \'e hiperel\ih tica. Vamos agora estudar outra fun\c c\~ao hiperel\ih tica em $T$, que \'e dada por $\cH(z)=-z+w_1+w_2$ (rota\c c\~ao em torno de $\ds\frac{w_1+w_2}{2}$). Seu conjunto de pontos fixos \'e exatamente $\ds\biggl\{\pm\frac{w_1+w_2}{2},$ $\ds\pm\frac{w_1-w_2}{2}\biggl\}$. Do Teorema 2.2, $\cH$ e $\wp$ induzem outra involu\c c\~ao em $\hC$, e como ambas s\~ao meromorfas, a involu\c c\~ao induzida tamb\'em ser\'a meromorfa. Do Teorema 2.1 esta involu\c c\~ao \'e uma transforma\c c\~ao de M\"obius. Como $\cH$ intercambia 0 e $w_1+w_2$ e deixa fixo $\ds\frac{w_1+w_2}{2}$, se $J$ \'e a induzida por $\wp$ temos:
\\

$J(0)=J(\wp(0))=\wp(\cH(0))=\wp(w_1+w_2)=\infty;$

$J(\infty)=J(\wp(w_1+w_2))=\wp(\cH(w_1+w_2))=\wp(0)=0;$

$J(i)=J(\wp(\ds\frac{w_1+w_2}{2}))=\wp(\cH(\frac{w_1+w_2}{2}))=\wp(\frac{w_1+w_2}{2})=i.$
\\

Assim, a transforma\c c\~ao de M\"obius associada \'e 
\BE
   \wp\to-\frac{1}{\wp}.
\EE

Uma vez que os pontos $\ds\pm\frac{w_1-w_2}{2}$ permanecem fixos pela $\cH$, neles temos $\wp^2=-1$. Como o grau da $\wp$ \'e 2 e $\cH(\ds\pm\frac{w_1+w_2}{2})\neq\cH(\pm\frac{w_1-w_2}{2})$, temos
\[
   \wp\biggl(\pm\frac{w_1-w_2}{2}\biggl)=-i.
\]

Al\'em disso, (3) implica que todo segmento em $T$ com extremos $p,q$, e cujo centro \'e um ponto fixo de $\cH$, satisfaz $\wp(q)=-1/\wp(p)$. Veja a figura abaixo.
\begin{figure}[h]
\center
\includegraphics[scale=0.8]{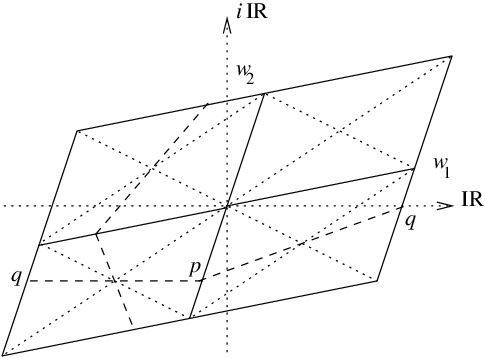}\\
\centerline{Figura 2.2: $p,q$ diagonalmente opostos no sub-reticulado, $\ds\wp(q)=-\frac{1}{\wp(p)}$.}
\end{figure}

De fato, como $J$ \'e induzida de $\cH$ pela $\wp$, temos $(\wp\circ\cH)(p)=(J\circ\wp)(p)=-1/\wp(p)$. Assim, basta mostrarmos que $\cH(p)=q$. Tomamos $v$ um ponto fixo de $\cH$ e $p,q$ os extremos de um segmento com ponto m\'edio $v$. Os pontos $p,q$ e $v$ satisfazem: $|\overrightarrow{pq}|=2|\overrightarrow{vq}|$, donde $p=2v-q$. Como $2v-w_1-w_2\in G$, temos o resultado para $v\in\biggl\{\ds\pm\frac{w_1+w_2}{2},\pm\frac{w_1-w_2}{2}\biggl\}$. 
\\

No caso especial de um toro retangular, a menos de biholomorfismo ou anti-biholomorfismo podemos assumir que $w_1\in\R_+$ e $w_2\in i\R_+$. Algumas involu\c c\~oes adicionais tamb\'em v\~ao satisfazer as hip\'oteses do Teorema 2.2, a saber:
\\
\\
\centerline{$I_1:z\to\ovl{z}$;}
\centerline{$I_2:z\to-\ovl{z}$;}
\centerline{$I_3:z\to\ovl{z}+w_2$;}
\centerline{$I_4:z\to-\ovl{z}+w_1$.}
\\

Elas est\~ao representadas na Figura 2.3:
\begin{figure}[h]
\center
\includegraphics[scale=0.8]{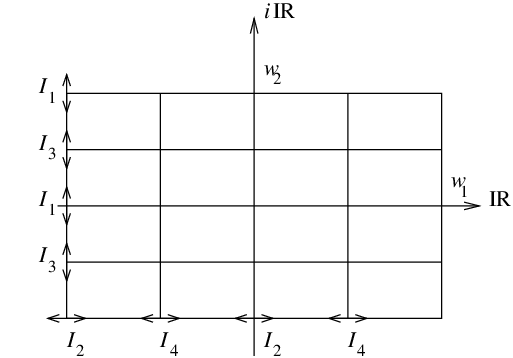}
\centerline{Figura 2.3: As involu\c c\~oes $I_{1,2,3,4}$ no toro retangular $T$.}
\end{figure}

Junto com $\wp$, elas induzem involu\c c\~oes anti-holomorfas $J_{1,2,3,4}$ em $\hC$ que identificaremos com o aux\ih lio do Teorema 2.2.
\\

Consideramos $I_1$ e $I_2.$ Elas fixam alguns pontos especiais como $0,w_1+w_2$ e intercambiam $\ds\frac{w_1+w_2}{2}$ e $\ds\frac{w_1-w_2}{2}.$ Uma vez que $\wp(0)=0$, $\wp(w_1+w_2)=\infty$ e $\wp(\ds\frac{w_1+w_2}{2})=-\wp(\frac{w_1-w_2}{2})=i$, elas induzem a mesma involu\c c\~ao em $\hC$, que fixa $0,\infty$ e intercambia $i$ com $-i$. Do Teorema 2.1, esta involu\c c\~ao \'e $\wp\to\ovl{\wp}$. Em particular, isto significa que a imagem de $\wp$ do conjunto dos pontos fixos de $I_{1,2}$ \'e {\it real}. Assim, $\wp(w_1)=\tan\af$ para algum $\af\in(-\pd,\pd)\setminus\{0\}$. Pela Figura 2.3, $w_1$ e $w_2$ s\~ao os extremos de um segmento com ponto m\'edio $\ds\frac{w_1+w_2}{2}$, e conclu\ih mos que $\wp(w_2)=-\cot\af$. Como $\tan\af=-\cot(\af+\pd)$, a menos de um anti-biholomorfismo podemos escolher o nosso toro de modo que $\af>0$. Dessa forma, consideramos $\af\in(0,\pd)$.
\\

A involu\c c\~ao $I_4$ intercambia $0$ com $w_1$ e $w_1+w_2$ com $w_2$. Logo, se $J_4$ for a involu\c c\~ao induzida por $\wp$ (de $I_4$), ela ser\'a uma transforma\c c\~ao de M\"obius anti-holomorfa que satisfaz:
\\

$J_4(0)=J_4(\wp(0))=\wp(I_4(0))=\wp(w_1)= \tan \af;$

$J_4(\infty)=J_4(\wp(w_1+w_2))=\wp(I_4(w_1+w_2))=\wp(w_2)= -\cot\af;$

$J_4(-\cot \af)=J_4(\wp(w_2))=\wp(I_4(w_2))=\wp(w_1+w_2)=\infty.$
\\

Portanto $J_4$ \'e dada por 
\[
   \wp\to\frac{\tan\af-\ovl{\wp}}{1+\tan\af\cdot\ovl{\wp}}.
\]

Os c\'alculos para determinar a involu\c c\~ao induzida por $I_3$ s\~ao an\'alogos. A pr\'oxima figura representa $\wp(T)$ com as suas involu\c c\~oes induzidas (foram dados os mesmos nomes).
\\
\begin{figure}[h]
\center
\includegraphics[scale=0.85]{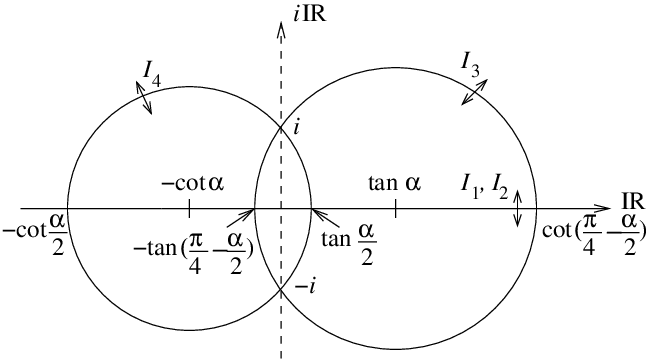}
\centerline{Figura 2.4: A imagem $\wp(T)$ com as involu\c c\~oes induzidas.}
\end{figure}

Usualmente escrevem-se os valores de $\wp$ sobre o toro, como na pr\'oxima figura.
\\
\begin{figure}
\center
\includegraphics[scale=1]{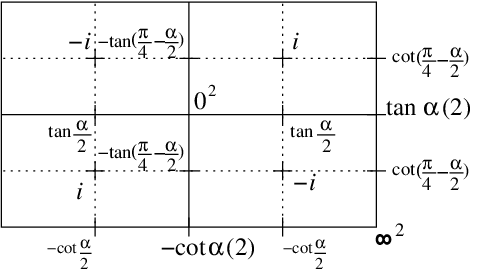}
\centerline{Figura 2.5: Os valores de $\wp$ no toro retangular $T$.}
\end{figure}

Agora estudaremos o caso especial do toro r\^ombico. A menos de um biholomorfismo ou anti-biholomorfismo podemos assumir que $Arg(w_1)\in(0,\pd)$ e $w_2=-\ovl{w}_1$. As involu\c c\~oes definidas anteriormente $I_{3,4}$ n\~ao s\~ao mais v\'alidas aqui, mas $I_{1,2}$ e outras duas involu\c c\~oes ainda v\~ao satisfazer as hip\'oteses do Teorema 2.1, a saber:
\\
\\
\centerline{$I_5:z\to  \overline{z}+w_1+w_2$;}
\centerline{$I_6:z\to -\overline{z}+w_1-w_2$.}
\\

Elas est\~ao representadas na Figura 2.6.
\\
\begin{figure}[h]
\center
\includegraphics[scale=0.8]{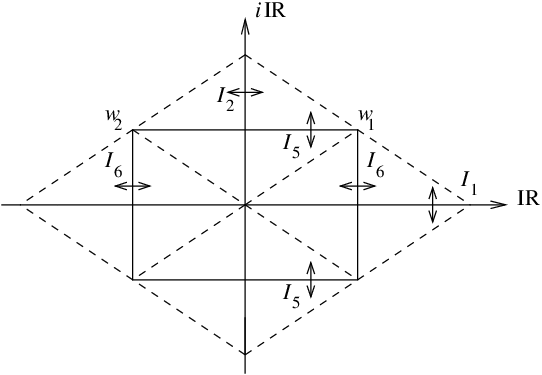}
\centerline{Figura 2.6: As involu\c c\~oes $I_{1,2,5,6}$ no toro r\^ombico $T$.}
\end{figure}

Em $\hC$, $I_{1,2}$ induzem a mesma involu\c c\~ao, a saber $\wp\to-\ovl{\wp}$. O mesmo ocorre com $I_{5,6}$ que tamb\'em induzem a mesma involu\c c\~ao: $\wp\to 1/\ovl{\wp}$. Uma conseq\"u\^encia importante da involu\c c\~ao $\wp\to-\ovl{\wp}$ em $\hC$ \'e que a imagem pela $\wp$ da diagonal de $T$ \'e o eixo imagin\'ario. Como escolhemos $\wp\biggl(\ds\frac{w_1+w_2}{2}\biggl)=i$, a imagem pela $\wp$ da diagonal horizontal de $T$ \'e $i\R_{-}$, enquanto que a imagem pela $\wp$ da diagonal vertical de $T$ \'e $i\R_+.$ Da involu\c c\~ao $\wp\to 1/\ovl{\wp}$ em $\hC$, temos que a imagem do sub-ret\^angulo representado na Figura 2.6 cobre $S^1$. Dessa forma, $\wp(w_1)=e^{i\rho}$ para algum $\rho\in(-\pd,\pd)$. Similarmente ao caso do toro ret\^angulo, sem perda de generalidade tomamos $\rho\ge 0$ e agora conside- ramos $\rho\in[0,\pd)$. Da Figura 2.6 conclu\ih mos que $\wp(w_2)=-e^{-i\rho}$.
\\

No caso particular do toro quadrado, que \'e simultaneamente r\^ombico e retangular, temos $\af=\pi/4$ e $\rho=0$. A pr\'oxima figura descreve os valores de $\wp$ diretamente sobre o toro.
\\
\begin{figure}[h]
\center
\includegraphics[scale=0.7]{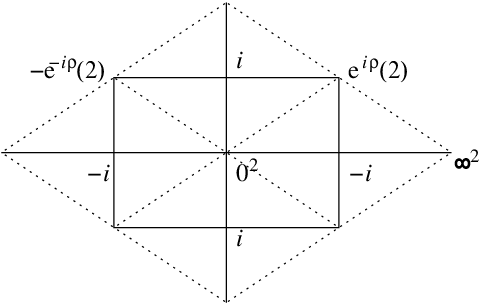}
\centerline{Figura 2.7: Os valores de $\wp$ no toro r\^ombico $T$.}
\end{figure}
\ \\
{\bf 3. A equa\c c\~ao alg\'ebrica do toro}
\\

O teorema seguinte \'e um dos resultados centrais da Teoria de Superf\ih cies de Riemman:
\\
\\
{\bf Teorema 2.3.} \it Toda superf\ih cie de Riemann compacta $X$ pode ser algebricamente descrita por duas fun\c c\~oes $f,g:X\to\hC$ que satisfazem uma equa\c c\~ao polinomial a duas indeterminadas $f$ e $g$, com coeficientes constantes em $\C$. Reciprocamente, toda equa\c c\~ao polinomial a duas indeterminadas $f$ e $g$, com coeficientes constantes em $\C$, admite uma superf\ih cie de Riemann compacta $X$ tal que $f,g$ s\~ao fun\c c\~oes de $X$ em $\hC$.\rm
\\

Com base neste teorema, vamos deduzir uma equa\c c\~ao alg\'ebrica para o toro $T=\C/G$. Na Figura 2.8(a), observe os \'unicos pontos de ramo da $\wp$, marcados com $\bullet$, onde $x=\wp(w_1)\in\C^*$. Como deg$(\wp)=2$, estes s\~ao os \'unicos pontos do toro em que $\wp$ assume tais valores. 
\\
\begin{figure}[h]
\center
\includegraphics[scale=0.7]{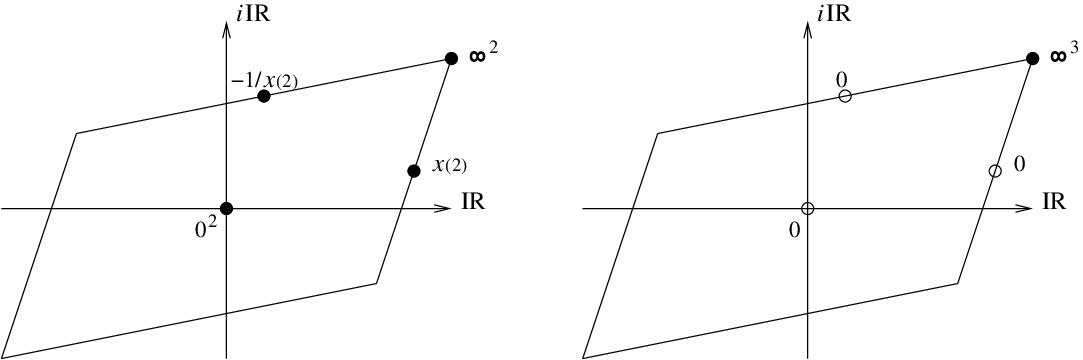}
\hspace{3.5cm}(a)\hspace{6.3cm}(b)
\centerline{Figura 2.8: (a) valores de $\wp$; (b) valores de $\wp'$.}
\end{figure}

Na Figura 2.8(b), marcamos apenas os p\'olos e zeros de $\wp'$. Note que esta possui um \'unico p\'olo (de ordem 3), sendo que seus zeros s\~ao todos simples, uma vez que deg$(\wp')=3$.
\\

Observe agora que a fun\c c\~ao $f=\wp(\wp-x)(\wp+1/x)$ tem exatamente os mesmos p\'olos e zeros que $g=\wp^{\prime 2}$. Numa vizinhan\c ca de cada cada ponto $\bullet$, considere os desenvolvimentos de Laurent para $f$ e $g$. Ent\~ao vemos que $f/g$ \'e uma fun\c c\~ao holomorfa de $T$ em $\C$, uma vez que n\~ao tem p\'olos. Pelo que j\'a discutimos na Se\c c\~ao 2.1, esta deve ser uma constante $c\in\C$. Al\'em disso, $c\ne 0$ pois $f\not\equiv 0$. Assim, obtemos uma equa\c c\~ao alg\'ebrica para $T$ dada por
\BE
   \wp'^2=c\wp(\wp-x)(\wp+1/x).
\EE

No caso particular do toro retangular, $x=\tan\af\in\R_+^*$ e (4) se reescreve como
\[
   \wp'^2=c\wp(\wp-\tan\af)(\wp+\cot\af).
\]

No caso particular do toro r\^ombico, $x=e^{i\rho}$ e (4) se reescreve como
\[
   \wp'^2=c\wp(\wp-e^{i\rho})(\wp + e^{-i\rho}).
\]
\ \\
{\bf 4. A constante $c$}
\\

Come\c camos por observar que quaisquer dois toros com mesmo quociente $\ds\frac{w_1}{w_2}$ s\~ao biholomorfos, e reciprocamente. Basta vermos que, se $\ds\frac{w_1}{w_2}=\frac{\tw_1}{\tw_2}$, ent\~ao existe $g:\C\to\C$ com $g(z)=az$ e $a\in\C^*$ tal que $g(w_1)=\tw_1$ e $g(w_2)=\tw_2$. Da Observa\c c\~ao do Teorema 2.2 temos que $T$ \'e biholomorfo a $\tT$. A rec\ih proca tamb\'em \'e v\'alida pois $\C$ \'e simplesmente conexo, e da an\'alise complexa temos que todo biholomorfismo de $\C$ em $\C$ \'e uma transforma\c c\~ao linear afim. Como podemos sempre tomar nossos reticulados pela origem, ela \'e da forma $z\mapsto az$ para algum $a\in\C^*$.
\\

Consideramos agora $T=\C/\G$ e $\tT=\C/\tilde{\G}$, onde $\tw_1=aw_1$ e $\tw_2=aw_2$ e $g(z)=az+b$ com $a,b\in\C$, $a\ne 0$, conforme a Figura 2.9. Sejam $\bt(t):=t(w_1+w_2)$, $0\le t\le 1$ e $\tilde{\bt}=t(\tw_1+\tw_2)+b=a\bt(t)+b$, $0\le t\le 1$. Ent\~ao $\wp(\bt(t))$ \'e uma curva que liga $0$ a $\infty$ e passa por $i=\wp(\bt(\m))$. Logo, $\wp(\bt(t))'\big|_\m=y\in\C^*$, sendo que $y$ depende apenas de $\ds\frac{w_1}{w_2}$. Por outro lado, 
\begin{eqnarray*}
   [\wp(\bt(t))'\big|_\m]^2=y^2 &\impl& \wp'(\bt(1/2))\cdot[\bt'(1/2)]^2=y^2 \\
   &\impl& c\,i(i-x)(i+1/x)\cdot(w_1+w_2)^2=y^2 \\
   &\impl& c=\frac{y^2/w_2^2}{(1+w_1/w_2)^2(x-1/x-2i)}.
\end{eqnarray*} 

\begin{figure}
\center
\includegraphics[scale=0.75]{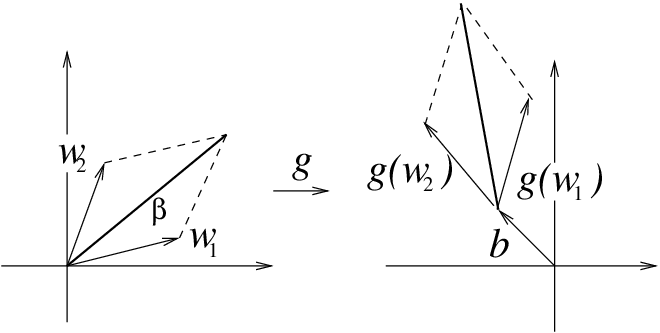}
\centerline{Figura 2.9: A aplica\c c\~ao $g$.}
\end{figure}

Assim, a rela\c c\~ao entre $c$ e o reticulado $\G$ para cada $u=\ds\frac{w_1}{w_2}$ fixo \'e \'unica, pois para cada par $(w_1,w_2),(\tw_1,\tw_2)$ com $\ds\frac{w_1}{w_2}=\frac{\tw_1}{\tw_2}=u$ temos: 
\[
   c|_{(w_1,w_2)}=c|_{(\tw_1,\tw_2)}\Longleftrightarrow(w_1+w_2)^2=(\tw_1+\tw_2)^2\Longleftrightarrow(w_1,w_2)=\pm(\tw_1,\tw_2),
\] 
que por sua vez define o mesmo reticulado $\G$. Logo, a menos de um biholomorfismo ou anti-biholomorfismo, todos os toros est\~ao representados num quadrante de c\ih rculo, conforme a Figura 2.10.
\eject
\begin{figure}[h]
\center
\includegraphics[scale=0.8]{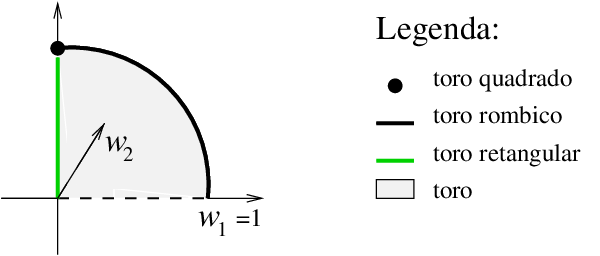}
\centerline{Figura 2.10: Representa\c c\~ao dos tipos de toros.}
\end{figure}
\eject
\ \\
{\huge\bf Cap\ih tulo 3}
\bigskip
\bigskip
\ \\
{\huge\bf A Fun\c c\~ao $\g$}
\bigskip
\bigskip
\bigskip

Apresentaremos, neste cap\ih tulo, um exemplo de fun\c c\~ao el\ih tica que pode ser obtido da fun\c c\~ao $\wp$ de Weierstrass Sim\'etrica atrav\'es de uma transforma\c c\~ao de M\"{o}bius adequada.
\\

Nosso exemplo \'e a fun\c c\~ao $\g$. Para constru\ih-la, considere a transforma\c c\~ao de M\"obius dada por $\ds Q(z)=e^{i\theta}\frac{i-z}{i+z}$, onde $\theta\in(0,\pi/2)$ ser\'a determinado como fun\c c\~ao de $\af$. Os principais valores-imagem de $\g=Q\circ\wp$ est\~ao descritos na Figura 3.1. 
\\

Os pontos de ramifica\c c\~ao de $\g$ s\~ao as pr\'e-imagens pela $\wp$ de $0,\tan\af,\infty$ e $-\cot\af$. Seus valores-imagem s\~ao, respectivamente, $e^{i \theta}, -e^{-i\theta},-e^{i \theta}$ e $e^{-i\theta}$.
\\
\begin{figure}[ht]
\centering
\includegraphics[scale=0.7]{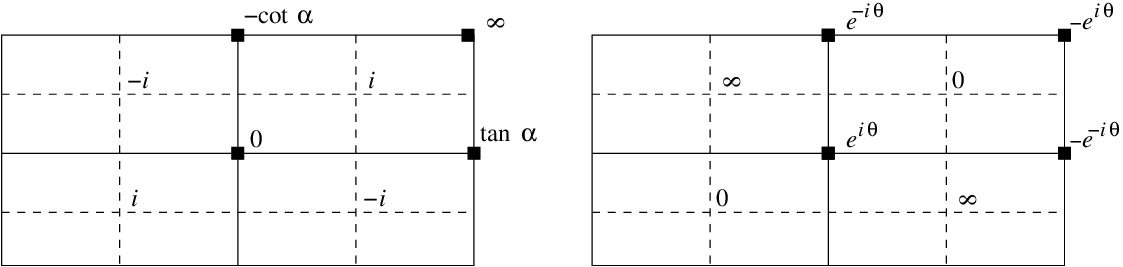}
\hspace{3.5cm}(a)\hspace{6.3cm}(b)
\centerline{Figura 3.1: a) valores de $\wp$; b) valores de $\g$.}
\end{figure}

Apresentamos, a seguir, um estudo detalhado da fun\c c\~ao $\g$ e de sua transforma\c c\~ao de M\"obius geradora. Neste estudo faremos uso de propriedades geom\'etricas do toro retangular para determinar os valores esbo\c cados na Figura 3.1(b).
\\

Seguindo procedimento an\'alogo ao explicado no Cap\ih tulo 2, podemos tamb\'em definir $\g$ de modo que:
\[
   \g(0)=e^{i\theta}, \ \ \g(w_1+w_2)=-e^{i\theta}, \ \ \g(\frac{w_1+w_2}{2})=0.
\]

Al\'em disso, 
\[
   \g(\pm\frac{w_1+w_2}{2})=Q(\wp(\pm \frac{w_1+w_2}{2}))=Q(i)=0
\]
e
\[  
   \g(\pm \frac{w_1-w_2}{2})=Q(\wp(\pm \frac{w_1-w_2}{2}))=Q(-i)=\infty.
\]

A reflex\~ao no eixo imagin\'ario $I_2(z)=-\bar{z}$ induz uma involu\c c\~ao anti-holomorfa de $\hC$ em $\hC$ que fixa $\g=e^{i \theta}$ e intercambia os valores $0$ e $\infty$. Logo, a involu\c c\~ao induzida $J_2$ \'e dada por  $\wp\to 1/\bar{\wp}$. Com isto, observamos que os pontos fixos de $J_2$ est\~ao sob o c\ih rculo unit\'ario $S^1$ (inclusive os pontos $\g(w_1),\g(w_2)$).
\\

A rota\c c\~ao de $180^\circ$ em torno de $\ds\frac{w_1+w_2}{2}$ \'e dada por $\cH(z)=-z +(w_1+w_2)$. A sua involu\c c\~ao induzida intercambia $e^{i\theta}$ com $e^{-i\theta}$, e fixa $0$. Assim, ela corresponde a $z\to-z$. Conseq\"uentemente, $\g(w_1)=-\g(w_2)$.
\\

Aplicando a involu\c c\~ao anti-holomorfa $I_3(z)=\bar{z}+w_2$, reflex\~ao na reta que passa por $z=\ds\frac{w_2}{2}$ e \'e paralela ao eixo real, observamos que a sua induzida $J_3$ fixa os pontos $\g=0$ e $\g=\infty$ e intercambia $e^{i\theta}$ com $\g(w_2)=Q(-\cot\af)$. Como $J_3(z)=\bar{z}$, temos $\g(w_2)=e^{-i\theta}$ e ent\~ao $\g(w_1)=-e^{-i\theta}$. Assim, obtemos a configura\c c\~ao esbo\c cada na Figura 3.1(b).
\\

Para estabelecermos uma rela\c c\~ao entre os \^angulos $\af$ e $\theta$, observamos que $Q(\tan\af)=-e^{-i\theta}$, donde

\begin{eqnarray*}
   e^{i\theta}\ds{\frac{i-\tan\af}{i+\tan\af}}=-e^{-i\theta} 
   &\Leftrightarrow & e^{i\theta}(i-\tan\af) =-e^{-i\theta}(i+\tan\af) \\
   &\Leftrightarrow & \tan\af\cdot(e^{i\theta}-e^{-i\theta})=i(e^{i\theta}+e^{-i\theta}) \\
   &\Leftrightarrow & \tan\af=\cot\theta \\
   &\Leftrightarrow & \af+\theta=\pi/2.
\end{eqnarray*}

Para determinarmos os valores de outros pontos do toro plano, observe que:
\begin{itemize}
\item $\ds\frac{w_2}{2}$ \'e fixado pela involu\c c\~ao $I_2$, logo $\ds\g(\frac{w_2}{2})\in S^1$;
\item $\ds\frac{w_2}{2}$ \'e fixado pela involu\c c\~ao $I_3$, logo $\ds\g(\frac{w_2}{2})$ \'e real.
\end{itemize}

Portanto, temos $\g(w_2/2)=\pm 1$. Uma vez que $\af\in(0,\pi/2)$, ent\~ao $\theta\in(0,\pi/2)$, donde $\g(w_2/2)=1$. Esse valor tamb\'em pode ser confirmado por $Q(-\tan(\pi/4-\af/2))=1$, pois
\[
   e^{i\theta}\cdot\frac{i+\tan(\pi/4-\af/2)}{i-\tan(\pi/4-\af/2)}=i\cdot e^{-i\af}\cdot\frac{i+\frac{\ds\cos\af}{\ds\sen\af+1}}{i-\frac{\ds\cos\af}{\ds\sen\af+1}}=i\cdot e^{-i\af}\cdot\frac{i+e^{i\af}}{i-e^{-i\af}}=1.
\]

Al\'em disso, considere as involu\c c\~oes anti-holomorfas $I_1(z)=\bar{z}$ (reflex\~ao no eixo real) e $I_4(z)= -\bar{z}+w_1$ (reflex\~ao na reta que passa por $z=\ds\frac{w_1}{2}$ e \'e paralela ao eixo imagin\'ario). Estas involu\c c\~oes induzem $J_1(z)=\ds\frac{1}{\bar{z}}$ e $J_4(z)=-\bar{z}$, respectivamente. Disto podemos inferir que:
\begin{itemize}
  \item $\ds\frac{w_1}{2}$ \'e fixado por $I_1$, logo $\g(\ds\frac{w_1}{2})\in S^1;$
  \item $\ds\frac{w_1}{2}$ \'e fixado por $I_4$, logo $\g(\ds\frac{w_1}{2})$ \'e imagin\'ario puro.
\end{itemize}

Portanto, temos $\g(\ds\frac{w_1}{2})=\pm i$. Uma vez que $\af\in(0,\pi/2)$, ent\~ao $\theta\in(0,\pi/2)$, donde $\g(w_1/2)=i$. Esse valor tamb\'em pode ser confirmado por $Q(\tan(\af/2))=i$, pois
\[
   e^{i\theta}\cdot\frac{i-\tan(\af/2)}{i+\tan(\af/2)}=i\cdot e^{-i\af}\cdot\frac{i-\frac{\ds\sen\af}{\ds\cos\af+1}}{i+\frac{\ds\sen\af}{\ds\cos\af+1}}=i\cdot e^{-i\af}\cdot\frac{1+e^{i\af}}{1+e^{-i\af}}=i.
\]

Devido \`a rela\c c\~ao de equival\^encia que define o toro, \'e poss\ih vel ``translad\'a-lo'', de forma que seus pontos de ramo estejam no interior do ret\^angulo, conforme a Figura 3.2.
\\

\begin{figure}[ht]
\centering
\includegraphics[scale=0.7]{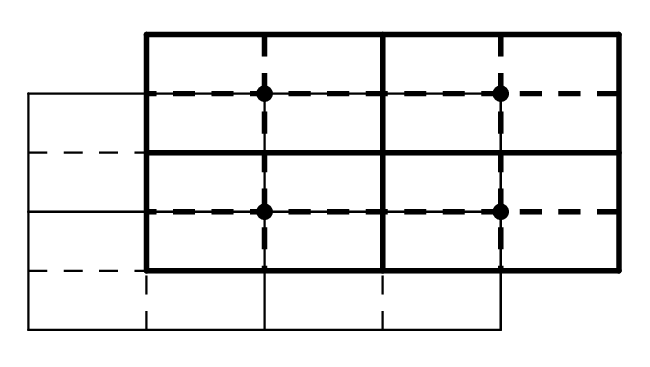}
\centerline{Figura 3.2: Diferentes maneiras de se representar o toro.}
\end{figure}

Dessa forma, \'e poss\ih vel representar os valores da fun\c c\~ao $\g: T \to \hC$ sob o toro conforme a Figura 3.3, descrita abaixo.
\\

\begin{figure}[ht]
\centering
\includegraphics[scale=0.9]{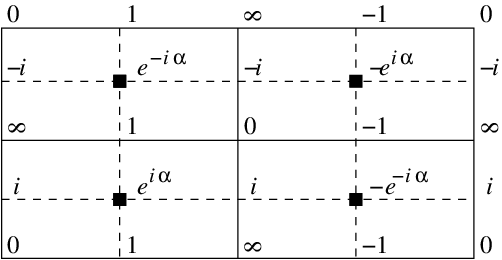}
\centerline{Figura 3.3: Valores da fun\c c\~ao $\g$.}
\end{figure}

Vamos agora determinar uma equa\c c\~ao alg\'ebrica para este toro. Para os pontos de ramo da fun\c c\~ao $\g$, precisamente $\g^{-1}(\{\pm e^{\pm i \theta}\})$, temos que express\~ao $\g^2+\g^{-2}-e^{2i\theta}-e^{-2i\theta}$ se anula. Considerando-se tamb\'em os p\'olos e zeros de $\g,\g'$ e $\g'/\g$, descritos na Figura 3.4, obtemos a equa\c c\~ao 
\[
   (\g'/\g)^2=c_0(\g^2+\g^{-2}-2\cos(2\theta)), \ c_0 \in \C^*.
\]
\begin{center}
\includegraphics[scale=0.48]{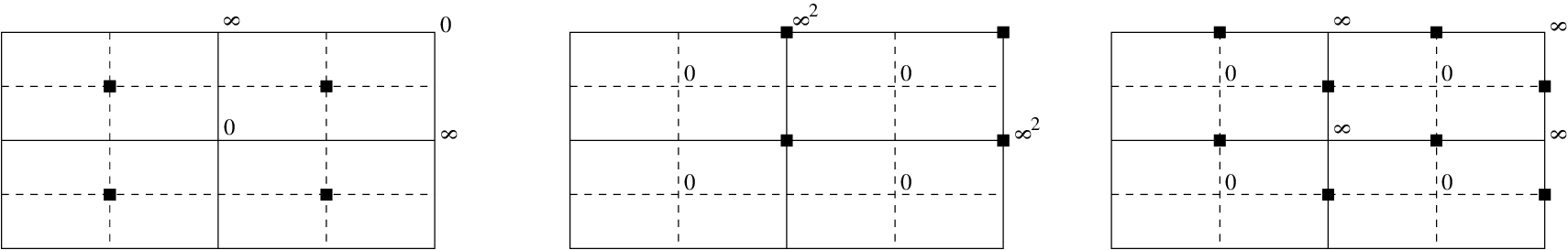}
a) \hspace{4.3cm} b)  \hspace{4.3cm} c)
\centerline{Figura 3.4: p\'olos e zeros de a) $\g$; b) $\g'$ e c) $\g'/\g$.}
\end{center}

Para uma curva $\delta(t)$ no toro que parametriza o eixo real do plano complexo quando aplicada a $\g,$ ou seja, $\g(\delta(t))=t\in\R$, teremos:
\[
   \g'(\delta(t))\cdot\delta'(t)=\ds\frac{d}{dt}(\g(\delta(t)))=\ds\frac{d}{dt}(t)=1.
\]

Portanto, $\g'(\delta(t))^2\delta'(t)^2=1$.
\\

Por outro lado, como $\delta$ \'e uma parametriza\c c\~ao do eixo real do toro, segue-se que $\delta'(t)=h(t)\in\R$. Deste modo, $\g'(\delta(t))^2=\ds\frac{1}{h(t)^2}\in\R_+$, donde  $\g'|_{\delta}\in\R$.
\\

Assim, $( \frac{\g'}{\g})^2 |_{\delta} \in \R_+$ e $(\g^2+\g^{-2} -2\cos 2\theta) |_{\delta} \in \R_+$, donde temos que $c_0$ \'e real positiva. Assim, $c_0$ determina apenas uma dilata\c c\~ao do reticulado, e fixamos $c_0=1$. Logo, uma equa\c c\~ao alg\'ebrica para $T$ \'e dada por 
\[
   \left(\frac{\g'}{\g}\right)^2=\g^2+\g^{-2}-2\cos 2\theta.
\]

Para finalizarmos este cap\ih tulo, apresentaremos outra rela\c c\~ao muito usa- da entre as fun\c c\~oes el\ih ticas $\wp$ e $\g$, mas que n\~ao pode ser obtida pela transforma\c c\~ao de M\"obius $Q$, pois agora transladamos o reticulado. Para obter a nova rela\c c\~ao, deduzimos da Figura 3.5 os p\'olos e zeros de 
\[
   \wp-\ds{\frac{1}{\wp}}+\tan\theta-\cot\theta=\wp-\ds{\frac{1}{\wp}}+\cot\af- \tan\af=(\wp+\cot\af)-(\ds{\frac{1}{\wp}}+\tan\af).
\]

\begin{center}
\includegraphics[scale=0.65]{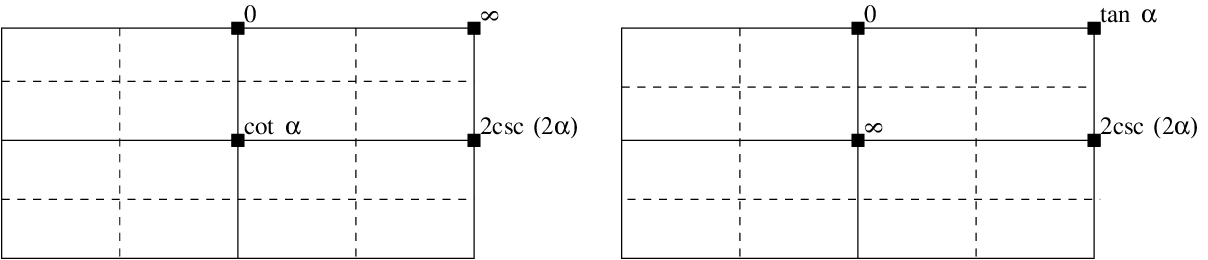}
\hspace{3.5cm}(a)\hspace{6.3cm}(b)
\centerline{Figura 3.5: P\'olos e zeros de a) $\cot\af+\wp$ e b) $\tan\af+1/\wp$.}
\end{center}

Logo, \'e imediato observar que os p\'olos e zeros de $(\wp-1/\wp+\tan\theta-\cot\theta)^{-1}$ s\~ao os mesmos (e com a mesma ordem) de $\g^2$ (veja a Figura 3.3). Portanto, estas equa\c c\~oes coincidem a menos de uma constante de proporcionalidade, isto \'e, 
\[
   \g^2=\ds\frac{c_1}{\wp-1/\wp+\tan\theta-\cot\theta} \ , \ c_1 \in \C^*.
\]

Mas para $\wp=\tan \frac{\af}{2}$ a fun\c c\~ao $\g$ assume o valor $i$, e como $\tan \frac{\af}{2}-\cot \frac{\af}{2}=$ $-2\cot \af$ temos que $c_1=\cot \af +\tan \af$, donde
\[
   \g^2=\ds\frac{\cot\af+\tan\af}{\wp-1/\wp+\tan\theta-\cot\theta}.
\]
\eject
\ \\
{\huge\bf Cap\ih tulo 4}
\bigskip
\bigskip
\ \\
{\huge\bf Aplica\c c\~oes}
\bigskip
\bigskip
\bigskip
\ \\
{\bf 1. Resultados cl\'assicos}
\\

Neste cap\ih tulo faremos uma aplica\c c\~ao da teoria de fun\c c\~oes el\ih ticas (re)cri- ando um exemplo de superf\ih cie m\ih nima. Come\c camos por introduzir alguns resultados sobre a teoria de superf\ih cies m\ih nimas. Para maiores detalhes vide [1], [2], [3] ou [4].
\\
\\
{\bf Teorema 4.1.} Seja $X:R\to\mathbb{E}$ uma imers\~ao isom\'etrica completa de uma superf\ih cie de Riemann $R$ sobre um espa\c co $\mathbb{E}$ ``flat'', completo e tridimensional. Se $X$ \'e m\ih nima e sua curvatura Gaussiana total $\int_{R} K dA$ \'e finita, ent\~ao $R$ \'e biholomorfa a uma superf\ih cie de Riemann compacta $\ovl{R}$ perfurada em um n\'umero finito de pontos $\{p_1,p_2,\cdots,p_s\}$.
\\
\\
{\bf Teorema 4.2.} (Representa\c c\~ao de Weierstra\ss). Seja $R$ uma superf\ih cie de Riemann, $g$ e $dh$ fun\c c\~ao e 1-forma meromorfas em $R,$ tais que os zeros de $dh$ coincidem com os zeros e p\'olos de $g$. Suponhamos que $X:R\to\mathbb{E}$, dada por 
\[
   X(p)=\int^p(\phi_1,\phi_2,\phi_3),\eh\eh(\phi_1,\phi_2,\phi_3)=\frac{1}{2}\biggl(\frac{1}{g}-g,\frac{i}{g}+ig,2\biggl)dh,
\]
esteja bem definida. Ent\~ao $X$ \'e uma imers\~ao m\ih nima conforme. Reciprocamente, toda imers\~ao m\ih nima conforme $X:R\to\mathbb{E}$ pode ser expressa como acima para alguma fun\c c\~ao $g$ e 1-forma $dh$ meromorfas.
\\
\\
{\bf Defini\c c\~ao 4.1.} O par $(g,dh)$ s\~ao os {\it dados de Weierstra\ss} e $\phi_{1,2,3}$ as {\it formas de Weierstra\ss} em $R$ da imers\~ao m\ih nima $X: R \to X(R) \subset \mathbb{E}.$
\\
\\
{\bf Teorema 4.3.} Nas hip\'oteses dos Teoremas 4.1 e 4.2, os dados de Weierstra\ss \  $(g,dh)$ se estendem meromorficamente sobre $\overline{R}.$
\\
\\
{\bf Defini\c c\~ao 4.2.} Um \it fim \rm de $R$ \'e a imagem de uma vizinhan\c ca perfurada $V_p$ de um ponto $p\in\{p_1,p_2,\cdots,p_s\}$ tal que $(\{p_1,p_2,\cdots,p_s\}\setminus p)\cap\ovl{V}_p=\emptyset$. O fim \'e mergulhado se sua imagem \'e mergulhada para uma vizinhan\c ca de $p$ suficientemente pequena.
\\
\\
{\bf Teorema 4.4.} (F\'ormula de Jorge-Meeks). Seja $X: R \to \mathbb{E}$ uma superf\ih cie m\ih nima regular completa de curvatura total finita $\int_RKdA$. Se os fins de $R$ s\~ao mergulhados, ent\~ao $$deg(g)=k+r-1$$ onde $k$ \'e o g\^enero de $\overline{R}=R\setminus\{p_1,p_2,\cdots,p_s\}$ e $r$ \'e o n\'umero de fins.
\\
\\
{\bf Exemplo 4.1.} O caten\'oide \'e um exemplo de superf\ih cie m\ih nima, cujos dados de Weierstra\ss \  s\~ao $g=z$, $dh=dz/z$, e a superf\ih cie de Riemann compacta \'e $\hC$. O caten\'oide possui dois fins, nos pontos $p_0=0$ e $p_1=\infty.$
\\

Um c\'alculo simples mostra que a imers\~ao m\ih nima do caten\'oide \'e dada por $\ds X(x,y)=\left(\frac{x}{x^2+y^2}+x,\frac{y}{x^2+y^2}+y,2\ln(\sqrt{x^2+y^2})\right)$.
\\

A fun\c c\~ao $g$ \'e a proje\c c\~ao estereogr\'afica da aplica\c c\~ao normal de Gau\ss \ $N:R\to S^2$ da imers\~ao m\ih nima $X$, ou seja, 
\[
   N=\frac{1}{|g|^2+1}(2 Re\{g\},2Im\{g\},|g|^2-1).
\]
Ela \'e um recobrimento (ramificado) de $\hC$. Al\'em disso, $\int_RKdA=-4\pi deg(g)$.
\\

O elemento de reta $ds$ de $X:R \to \mathbb{E}$ \'e dado por
\[
   ds=\frac{1}{2}\left(|g|+\frac{1}{|g|}\right)|dh|
\]
e a curvatura Gaussiana \'e dada pela seguinte f\'ormula:
\[
   K=-\biggl(\frac{2}{|g|+1/|g|}\biggl)^4\biggl|\frac{dg/g}{dh}\biggl|^2.
\]
\ \\
{\bf Teorema 4.5.} Se $\sigma$ \'e uma curva em $X(R)$ ent\~ao vale:
\begin{enumerate}
\itemsep = 0.0 pc
\parsep  = 0.0 pc
\parskip = 0.0 pc
\item [i)] $\sigma$ \'e assint\'otica se, e somente se, $\frac{dg}{g}(\sigma')\cdot dh(\sigma')\in i\R$;
\item[ii)] $\sigma$ \'e principal se, e somente se, $\frac{dg}{g}(\sigma')\cdot dh(\sigma')\in\R$.
\end{enumerate}
\ \\
{\bf Teorema 4.6.} Seja $\sigma$ uma curva anal\ih tica em $R$ tal que $g(\sigma)$ est\'a contida em um meridiano $e^{i\theta}\R$ ou no equador $S^1$, e $dh(\sigma')\subset\R$ ou $dh(\sigma')\subset i\R$. Ent\~ao $\sigma$ \'e uma geod\'esica de $R$ e uma linha de simetria.
\\
\\
{\bf 2. Exemplo de superf\ih cie m\ih nima}
\\

Considere a superf\ih cie da Figura 4.1, que chamaremos {\it campo de cate\'oides}. Trata-se de uma fam\ih lia de superf\ih cies m\ih nimas duplamente peri\'odicas, cujos membros denotaremos por $M$. 
\\
\begin{figure}[ht]
\centering
\includegraphics[scale=0.8]{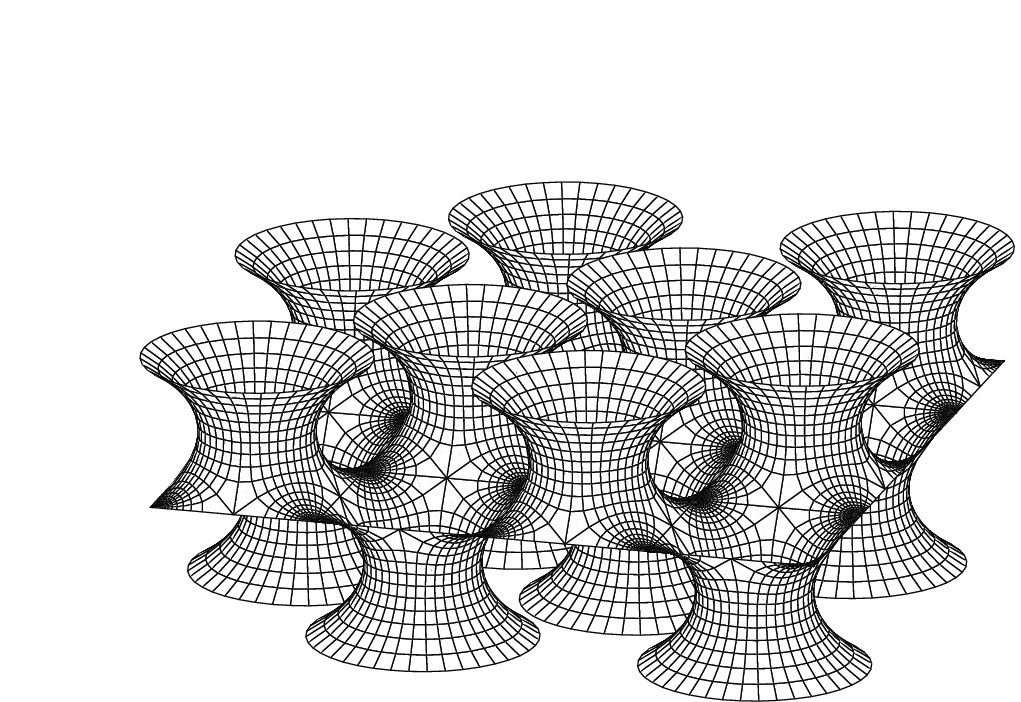}
\centerline{Figura 4.1: Um campo de caten\'oides.}
\end{figure}

Apresentamos a seguir a constru\c c\~ao de $M$. Para a obten\c c\~ao desta superf\ih cie, utilizamos o {\it m\'etodo da constru\c c\~ao reversa} desenvolvido por Hermann Karcher na d\'ecada de 80 (vide [1]). Seu m\'etodo consiste dos seguintes passos:
\begin{enumerate}
\itemsep = 0.0 pc
\parsep  = 0.0 pc
\parskip = 0.0 pc
 \item Esbo\c co da Superf\ih cie M\ih nima $M$;
 \item Compactifica\c c\~ao $\overline{M}$ de $M;$
 \item Hip\'oteses de Simetria;
 \item Equa\c c\~ao Alg\'ebrica de $\overline{M};$ 
 \item Obten\c c\~ao dos Dados de Weierstra\ss \  de $\overline{M};$
 \item Verifica\c c\~ao das Involu\c c\~oes de $\overline{M}$ e das Simetrias das Hip\'oteses;
 \item An\'alise dos Per\ih odos;
 \item Verifica\c c\~ao do Mergulho da Pe\c ca Fundamental.
\end{enumerate}

Consideramos uma superf\ih cie topol\'ogica que admita uma estrutura Riemanniana e seja isometricamente imersa em $\R^3$ com as seguintes propriedades: 1) a superf\ih cie imersa \'e m\ih nima; 2) ela \'e gerada pela rota\c c\~ao em torno de um segmento de reta junto com um grupo de transla\c c\~ao horizontal, ambos aplicados sobre uma pe\c ca fundamental $P$, onde $P$ \'e uma superf\ih cie com bordo e um fim catenoidal; 3) A proje\c c\~ao de $\deh P$ ortogonalmente sobre $x_3=0$ consiste de um quadrado $Q$ conforme a Figura 4.2(a).
\\

Consideramos os eixos $Ox_1$ e $Ox_2$ conforme a Figura 4.2(b) e definimos $x_3=x_1\wedge x_2$. Podemos interpretar $P$ como sendo a metade de um caten\'oide apoiado sobre um quadrado. Observa-se que, para a superf\ih cie quocientada pelo seu grupo de transla\c c\~ao, o seu topo catenoidal e a sua base catenoidal n\~ao se encontram no mesmo eixo. Esta superf\ih cie duplamente peri\'odica possui auto-intersec\c c\~oes, mas elas ocorrem somente nos fins catenoidais.
\\
\begin{figure}[ht]
\centering
\includegraphics[scale=0.6]{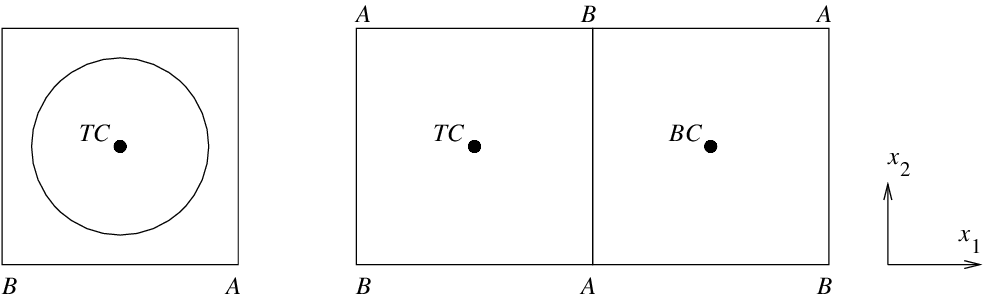}

\hspace{-3.0cm}(a)\hspace{4.2cm}(b) 

\centerline{Figura 4.2: a) o quadrado $Q$; b) a superf\ih cie compactificada $\overline{M}.$}
\end{figure}

Considerando o quociente de $M$ pelo seu grupo de transla\c c\~oes, seguido da compactifica\c c\~ao dos seus fins catenoidais ($TC$ - caten\'oide topo e $BC$ - caten\'oide base), obtemos uma superf\ih cie topol\'ogica compacta, que denotamos $\ovl{M}$ (vide Figura 4.2(b)). Desta figura vemos facilmente que o g\^enero de $\ovl{M}$ \'e 1, ou seja, $\ovl{M}$ \'e um toro.
\\

At\'e aqui, completamos os itens 1 e 2 do m\'etodo de Karcher. Para o item 3, vamos assumir que os trechos $A\to B$ formam as retas da superf\ih cie, enquanto que os caminhos $B\to TC$ e $A\to BC$ est\~ao contidos em planos paralelos ao plano $x_1=x_2$, e os caminhos $B\to BC$ e $A\to TC$ est\~ao em planos paralelos ao $x_1=-x_2.$
\\

Uma vez que as rota\c c\~oes de $180^\circ$ em torno dos trechos $A\leftrightarrow B$ fixam cada qual apenas uma componente conexa de $\ovl{M}$, este toro \'e r\^ombico. Note que o correspondente reticulado $\G$ de $\ovl{M}$ {\it n\~ao \'e paralelo} aos eixos $Ox_{1,2}$, mas sim rotacionado de $45^\circ$ em rela\c c\~ao a estes eixos. Al\'em disso, as reflex\~oes nos caminhos $B\to TC\to B$ e $A\to BC\to A$ fixam estas duas componentes conexas, donde temos que o toro tamb\'em \'e retangular. Conseq\"uentemente, $\ovl{M}$ \'e um toro quadrado.
\\

\begin{center}
 \includegraphics[scale=0.7]{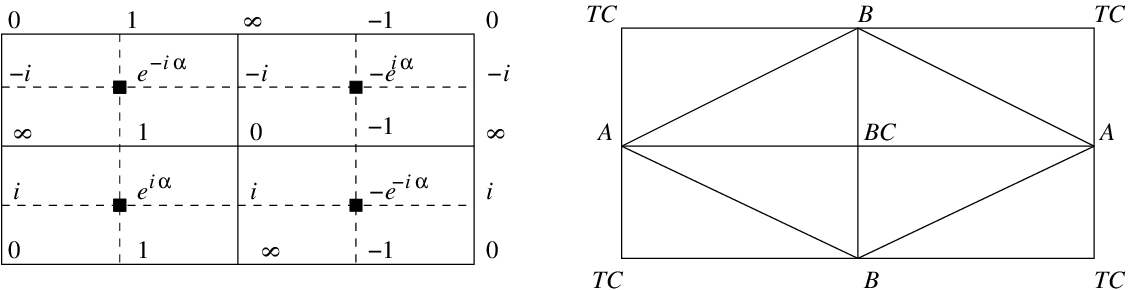}

\hspace{-0.5cm}(a) \hspace{6.5cm}(b) 

\centering{Figura 4.3: a) Valores de $\g$; b) correspond\^encia de pontos de $M$ no toro.}
\end{center}

Argumentos anal\ih ticos tamb\'em ir\~ao mostrar que trata-se de um toro quadrado.
\\

Do Cap\ih tulo 3 vimos que existe uma fun\c c\~ao $\g$ definida no toro, que assume os valores descritos na Figura 4.3(a). Seus pontos de ramo s\~ao indicados por $\blacksquare$, e $\af\in(0,\pi/2)$. As retas pontilhadas cobrem $S^1\subset\C$. Observe que o toro \'e quadrado se, e somente se, $\af=\pi/4$. Uma equa\c c\~ao alg\'ebrica para este toro \'e dada por 
\[
   \left(\frac{\g'}{\g}\right)^2=c_0 \left(\g^2+\g^{-2} -2\cos 2\alpha \right), \ c_0 >0.
\]

Tomando uma curva $\bt(t)$ no toro que parametriza a diagonal do plano complexo quando aplicada a $\g$, ou seja, $\g(\bt(t)) = e^{i\pi/4}t$, $t^{-1}\in (-1,1)$, temos
\[
   \g'(\bt(t))\cdot\bt'(t)=\frac{d}{dt}(\g(\bt(t)))=\frac{d}{dt}(e^{i\pi/4}t)=e^{i\pi/4}.
\]
Portanto, $\g'(\bt(t))^2\bt'(t)^2=i$.
\\

Observe a pr\'oxima figura.
\begin{center}
\includegraphics[scale=0.6]{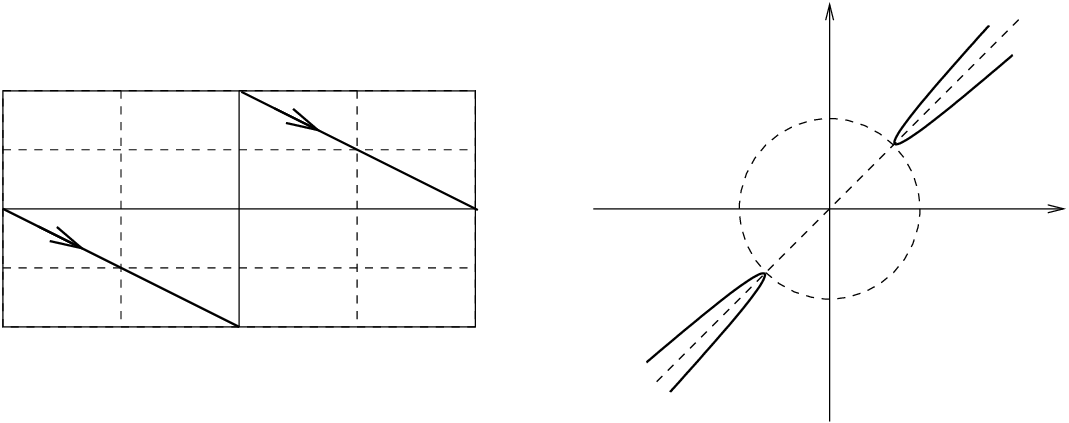}

\hspace{0.5cm} a) \hspace{5cm} b)
\centerline{Figura 4.4: a) Curva $\bt$ sobre o toro; b) imagem de $\bt$ por $\g$.}
\end{center}

Como $\bt(t)$ \'e uma parametriaza\c c\~ao de uma diagonal do toro, temos que $\bt'(t) = f(t) e^{i\pi/4},$ onde $f(t)$ \'e uma fun\c c\~ao real. Deste modo:
\begin{eqnarray*}
\g'^2(\bt(t))\cdot\bt^{\prime 2}(t)=i &\see & \g'^2(\bt(t))\cdot f^2(t)\cdot(-i)=i\\
& \see & \g'^2(\bt(t))=\frac{-1}{f^2(t)}\in\R_{-}\\
& \see & \g'(\bt(t))\in i\R.
\end{eqnarray*} 

Assim, $(\frac{\g'}{\g})^2|_{\bt(t)}\in i\R_+$. Mas $(\g^2+\g^{-2} -2\cos 2\alpha) |_{\bt(t)} \in i\R_+$ se, e somente se, $ \alpha=\pi/4$. Dessa forma, o toro \'e quadrado.
\\

At\'e aqui, completamos os itens 3 e 4 do m\'etodo de Karcher. Vamos agora partir para o item 5, que \'e a obten\c c\~ao dos dados de Weierstra\ss \ de $\overline{M}$. Como a aplica\c c\~ao $g$ \'e a proje\c c\~ao estereogr\'afica da aplica\c c\~ao de Gau\ss \ em $M$, observando a Figura 4.1 vemos que $M$ possui vetores normais unit\'arios verticais em $A,B,TC$ e $BC$. Se fixarmos $g(\{BC,TC\})=0$ teremos $g(\{ A,B\})=\infty$. Como o g\^enero de $\ovl{M}$ \'e 1 e o n\'umero de fins do quociente de $M$ pelo seu grupo de transla\c c\~ao \'e 2, pela f\'ormula de Jorge-Meeks temos 
\[
   deg(g)=k+s-1=1+2-1=2,
\]
e estes s\~ao os \'unicos p\'olos e zeros da $g$, todos simples.
\\

Como $g$ e $\g$ possuem exatamente os mesmos p\'olos e zeros, ent\~ao $g=c_1\g$ com $c_1\in\C^*$.
\\

No caminho $TC \to B$  a imagem da $g$ est\'a contida no plano $x_1=x_2,$ logo $g(\{TC \to B\}) \in e^{i\frac{\pi}{4}} \R_+,$ enquanto que $\g$ \'e real positiva neste caminho. Assim, uma equa\c c\~ao alg\'ebrica para $g$ \'e dada por 
\[
   g=c\,e^{i\pi/4}\,\g, \ c>0.
\]

Para determinarmos uma equa\c c\~ao da diferencial $dh$ observamos que, do Teorema 4.2, se $M$ \'e a imers\~ao m\ih nima completa de $\ovl{M}\setminus\{BC,TC\}$, ent\~ao todos os zeros e p\'olos de $dh$ s\~ao determinados pelos p\'olos e zeros da $g$ junto com os pontos regulares da $M$. Com isso, observamos que $dh(\{A,B\})=0$ e $dh(\{BC,TC\})=\infty$. Como 
\[
   -\chi(\ovl{M})=deg(dh)={\rm no. \ zeros}-{\rm no. \ p\acute{o}los},
\] 
a an\'alise est\'a conclu\ih da. Agora, $d\g/\g'$ representa uma forma holomorfa no toro. Ent\~ao, da Figura 4.5 obtemos 
\[
   dh=c_2\cdot\frac{1}{\g}\cdot\frac{d\g}{\g'}, \ c_2 \in \C^*.
\]

Queremos que a imers\~ao m\ih nima de $\bt(t)$ em $\R^3$ seja uma reta horizontal. Deste modo, devemos ter $Re\int_{\bt(t)} dh = 0$, ou seja, $dh|_{\bt(t)} \in i\R$. Assim, escolhemos a constante de propor\c c\~ao igual a 1, pois qualquer outro real s\'o reescalonaria a superf\ih cie, com ou sem invers\~ao pela origem. Logo, 
\[
   dh=\frac{1}{\g}\cdot\frac{d\g}{\g'}.
\]

\begin{center}
\includegraphics[scale=0.7]{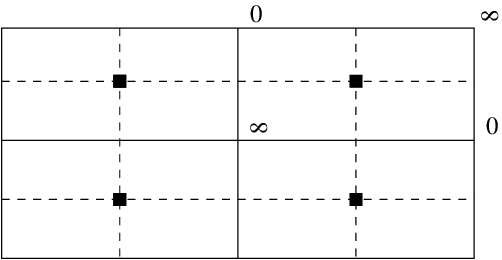}
\centerline{Figura 4.5: Zeros e pol\'os de $1/\g$ e $dh$.}
\end{center}

De posse destas informa\c c\~oes, definimos $M$ como a imers\~ao m\ih nima completa em $\R^3$ de $\overline{M} \setminus  \g^{-1}(\{0\}),$ sendo que a equa\c c\~ao alg\'ebrica de $M$ \'e dada por $(\g'/\g)^2=\g^2+1/\g^2$, e os dados de Weierstra\ss \  da superf\ih cie est\~ao definidos por $g=c\,e^{i\pi/4}\,\g$, $c>0$, e $dh=\ds\frac{d\g}{\g\g'}.$
\\

Estamos agora no item 6 do m\'etodo de Karcher, a \it verifica\c c\~ao das hip\'ote- ses de simetria. \rm Para a curva $\bt(t)$ vimos que $dh|_{\bt(t)} \in i \R$ e $g|_{\bt(t)} \in i\R$. Do Teorema 4.6, temos que $\bt$ \'e uma geod\'esica e uma linha de simetria. Do Teorema 4.5, a imers\~ao de $\bt$ representa uma reta na superf\ih cie, pois $\ds\frac{dg}{g}\,dh|_{\bt(t)} \in i \R.$
\\

Para o caminho $A \to TC$ tomemos $\rho(t)=it, t \in \R.$ Facilmente se verifica $\g'|_{\rho(t)} \in \R$ e al\'em disso, $dh|_{\rho(t)} \in \R$ e $g|_{\rho(t)} = i\,c\,e^{i \pi/4}t$, donde temos um meri- diano. Assim, $\rho$ \'e uma geod\'esica e linha de simetria. Como $\ds\frac{dg}{g}\,dh|_{\rho(t)} \in \R$, a imers\~ao \'e uma curva plana.
\\

Por argumentos an\'alogos pode-se provar que todas as demais simetrias se verificam.
\\

Para toda curva fechada em $\overline{M} \setminus \{BC, TC\}$, vamos analisar o vetor per\ih odo $Re\oint\phi_{1,2,3}$, item 7 do m\'etodo de Karcher.
\\

Tomando uma curva fechada simples em $\overline{M} \setminus \{BC, TC\}$ em torno de $TC,$ esta \'e homot\'opica a uma curva sim\'etrica em rela\c c\~ao aos trechos $A \to TC$ e $B \to TC$. Portanto, deve ser ortogonal a ambos os planos que cont\^em estas curvas. Como estes planos s\~ao perpendiculares, temos que $Re \oint \phi_{1,2,3} = 0$ sobre esta curva. De forma an\'aloga \'e poss\ih vel mostrar que $Re \oint \phi_{1,2,3} = 0$ sobre curvas fechadas simples em $\overline{M} \setminus \{BC, TC\}$ em torno de $BC.$
\\

\begin{center}
 \includegraphics[scale=0.7]{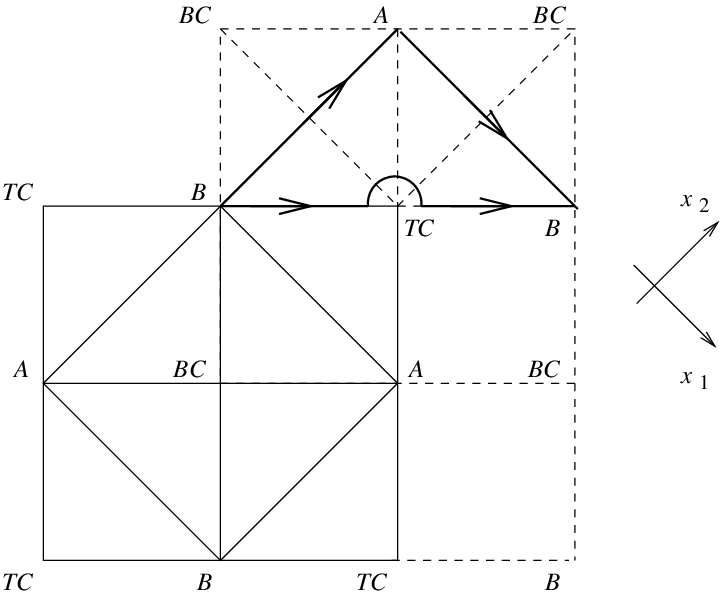}

\centering{Figura 4.6: Per\ih odo ao longo da curva $\sigma$.}
\end{center}

Considere a curva $\sigma$ homot\'opica a $B \to TC \to B$, com desvio \ em \ $TC$, descrita na Figura 4.6. Observe que o trecho $B \to A$ \'e sim\'etrico pela geod\'esica planar $BC \to TC,$ donde $Re \int_{B \to A} \phi_{1,3}=0.$ Por outro lado,
\[
   Re\int_{B \to A}\phi_2 =2Re\int_1^\infty i(g^{-1}+g)dh|_{\g(t)=te^{-i\frac{\pi}{4}}}=\pm 2\int_1^\infty\frac{(ct+(ct)^{-1})dt}{t(t^4-1)^{1/2}}\ne  0,
\]
o que nos d\'a um per\ih odo na dire\c c\~ao de $Ox_2$.
\\

Como provamos anteriormente que $A \to TC$ \'e uma geod\'esica de reflex\~ao, temos $|Re\int_{A \to B} \phi_1| = |Re\int_{B \to A}\phi_2| = \ld$. Portanto, $|Re\int_{B \to TC \to B} \phi_{1,2}|=\ld \sqrt{2} \neq 0.$ Por causa das simetrias que provamos existir, tamb\'em temos $|Re\int_{B \to BC \to B} \phi_{1,2}|=\ld \sqrt{2}.$
\\

Dessa forma, a superf\ih cie \'e duplamente peri\'odica.
\\

Finalmente, vamos mostrar que a pe\c ca fundamental da imers\~ao \'e mergulhada, \'ultimo item do m\'etodo de Karcher.
\\

O objetivo \'e verificar que as auto-interse\c c\~oes de $M$ ocorrem somente nos fins catenoidais.

\begin{center}
 \includegraphics[scale=0.7]{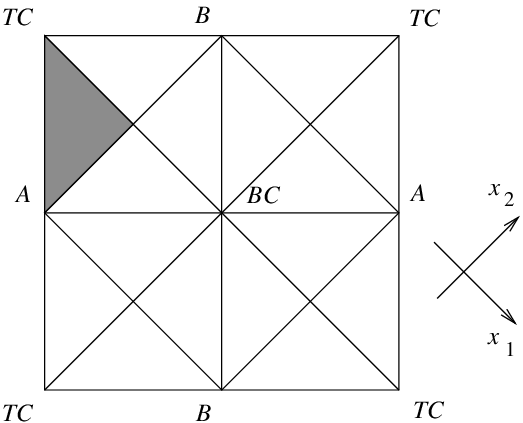}

\centering{Figura 4.7: Dom\ih nio $\cR$.}
\end{center}

Considere a regi\~ao $\cR$ descrita na Figura  4.7. A imagem $\mathcal{D}$ desta regi\~ao pela aplica\c c\~ao $g$ \'e um setor do plano complexo descrito na Figura 4.8(a). Como $g$ \'e a proje\c c\~ao estereogr\'afica da normal de Gau\ss, temos que $N(\cR)$ est\'a totalmente contida em um dos hemisf\'erios de $S^2$. Da imagem $\mathcal{D}$, conclu\ih mos que a restri\c c\~ao de $(x_2,x_3)$ ao interior de $\cR$ \'e uma imers\~ao, cuja imagem est\'a descrita na Figura 4.8(b).
\begin{center}
\includegraphics[scale=0.6]{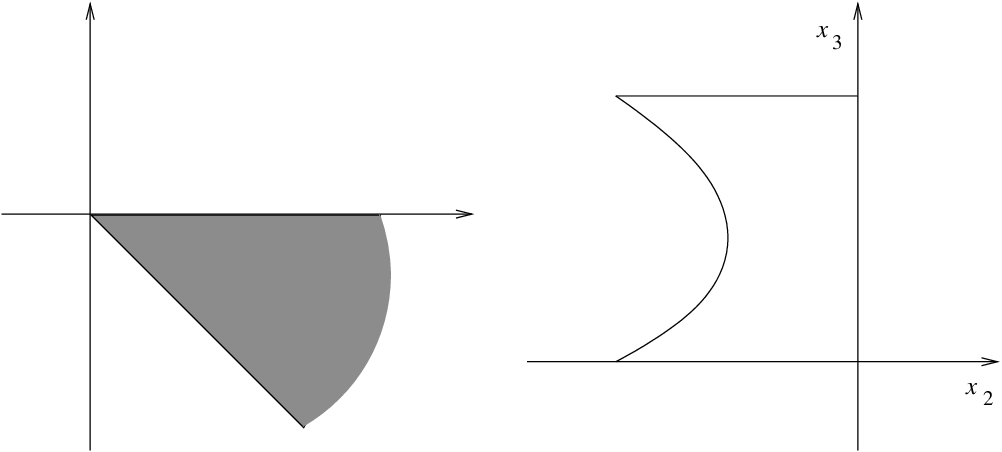}
\hspace{3.5cm}(a)\hspace{6.3cm}(b)
\centerline{Figura 4.8: a) imagem $\mathcal{D}=g(\cR)$; b) proje\c c\~ao $(x_2,x_3)|_{\cR}$.}
\end{center}

\begin{center}
\includegraphics[scale=0.9]{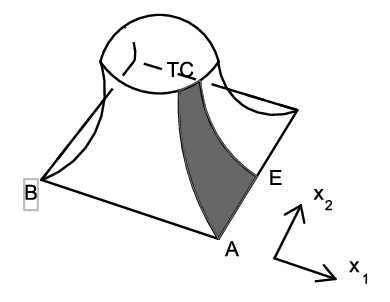}
\includegraphics[scale=0.9]{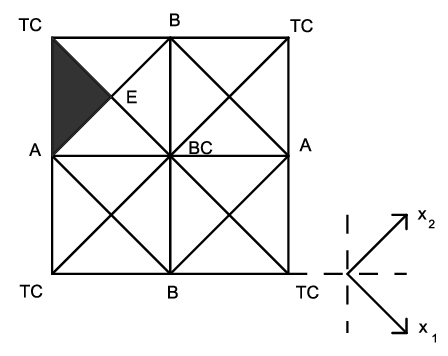}
\hspace{3.5cm}(a)\hspace{6.3cm}(b)
\centerline{Figura 4.9:a) dom\ih nio fundamental de $P$; b) dom\ih nio fundamental no toro.}
\end{center}

Com mencionamos anteriormente, $P$ \'e a pe\c ca fundamental de $M$. A regi\~ao hachurada na Figura 4.9(a) representa a imagem pela imers\~ao m\ih nima do dom\ih nio fundamental $\cR$, descrito pela Figura 4.9(b).
\\

Observe que o toro \'e dividido em 16 regi\~oes que s\~ao congruentes a $\cR$ e, por sua vez, a pe\c ca fundamental tamb\'em possui 8 regi\~oes que s\~ao obtidas da regi\~ao hachurada (Figura 4.9(a)) atrav\'es de simetrias da superf\ih cie. Como o quociente de $M$ pelo seu grupo de transla\c c\~ao \'e gerado por $P$ e a sua imagem por uma rota\c c\~ao de $180^{\circ}$, tamb\'em obtemos 16 subregi\~oes que formam este quociente.
\\

Como $g$ \'e uma aplica\c c\~ao aberta dada por $g=c\,e^{i\pi/4}\,\g$, $c>0$, $g$ aplica $\cR$ na regi\~ao descrita pela Figura 4.8(a). De fato, seja $\nu_1(t)$ uma curva no toro que descreve o caminho $TC\to A.$ Logo, $\g(\nu_1(t))=-it$, $t\in\R_+$. Aplicando esta curva em $g$ obtemos 
\[
   g(\nu_1(t))=c\,e^{i\frac{\pi}{4}}\,\g(\nu_1(t))=c\,e^{i\frac{\pi}{4}}\,\cdot(-it)=c\,e^{-i\frac{\pi}{4}}\,t.
\]

Tomando $\nu_2(t)$ uma curva no toro que descreve o caminho $TC \to E$, de modo que $\g(\nu_2(t))=e^{-i\frac{\pi}{4}}t$, $t\in [0,1]$, obtemos 
\[
   g(\nu_2(t))=c\,e^{i\frac{\pi}{4}}\,\g(\nu_2(t))=c\,e^{i\frac{\pi}{4}}\,e^{-i\frac{\pi}{4}}t=c\,t.
\]

Analogamente, se $\nu_3(t)$ descreve o caminho $E\to A$, ent\~ao $\g(\nu_3(t))=e^{-i\frac{\pi}{4}}\,t$, $t\in [1,\infty]$, donde temos
\[
   g(\nu_3(t))=c\,e^{i\frac{\pi}{4}}\,\g(\nu_3(t))=c\,e^{i\frac{\pi}{4}}\,e^{-i\frac{\pi}{4}}t=c\,t.
\]

Como $g$ \'e aberta, mas tem grau 2, ent\~ao $g(\cR)$ s\'o pode ser a regi\~ao indicada na Figura 4.8(a). Assim, $N(\cR)$ est\'a contido num hemisf\'erio, ou seja, existe uma dire\c c\~ao em que a proje\c c\~ao ortogonal do dom\ih nio fundamental \'e uma imers\~ao, e por conveni\^encia tomamos na dire\c c\~ao de $Ox_1$. Dessa forma, $(x_2,x_3):\cR\to\R^2$ \'e uma imers\~ao quando restrita ao interior de $\cR$.
\\

A proje\c c\~ao do dom\ih nio fundamental no plano $x_1=0$ possui duas alternativas que est\~ao indicadas na Figura 4.10.
\\

\begin{center}
\includegraphics[scale=0.9]{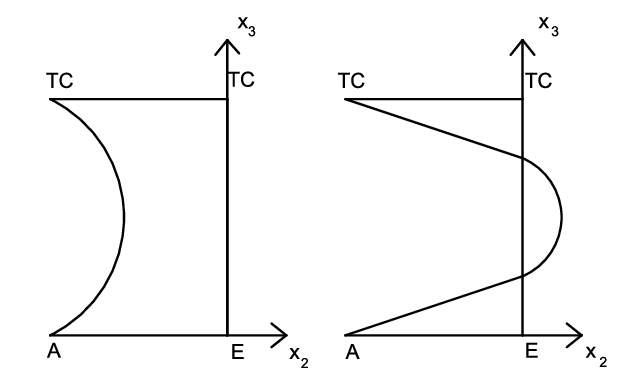}

\hspace{1.0cm}(a)\hspace{4.0cm}(b)
\centerline{Figura 4.10: poss\ih veis proje\c c\~oes do dom\ih nio fundamental em $x_1=0$.}
\end{center}

A curva $TC \to A$ \'e convexa. De fato, apresentamos na Se\c c\~ao 1 deste cap\ih tulo a f\'ormula da curvatura gaussiana para superf\ih cies m\ih nimas. Assim, $K=0$ se, e somente se, $dg$ se anula em $TC\to A$. Veremos que isto n\~ao ocorre na curva $TC\to A$.
\\

Basta observarmos que $deg(dg)=-\chi_{\bar{M}}=0$, e que $dg$ possui dois p\'olos com multiplicidade dois, logo o n\'umero de zeros de $dg$ deve ser quatro.
\\

Tamb\'em, para o caminho $TC \to BC$, al\'em de $g(\{TC, BC\})=0$ temos $g(TC\to BC)\in e^{-\frac{i\pi}{4}}\R_+$. Logo existe um ponto em $TC \to BC$ com $dg=0.$ Pelas simetrias da superf\ih cie, este ponto \'e aplicado sobre outros tr\^es pontos que tamb\'em s\~ao zeros da $dg.$ Portanto, a curvatura em $TC \to A$ nunca se anula, sendo esta curva convexa.
\\

Como $TC \to A$ \'e uma geod\'esica planar e, portanto, uma linha de curvatura, os contornos da Figura 4.10 est\~ao coerentes. Uma vez que $(x_2,x_3)$ \'e uma imers\~ao, ela \'e uma aplica\c c\~ao aberta e cont\ih nua, logo s\'o pode realizar o contorno da Figura 4.10(a).
\\

Dessa forma, $(x_2,x_3)$ \'e uma imers\~ao cujo contorno \'e uma {\it curva mon\'otona}, ou seja, uma curva $C^1$ por partes tal que seu vetor tangente se anula apenas para um conjunto discreto de pontos. Al\'em disso, a regi\~ao \'e simplesmente conexa. Assim, $(x_1,x_2, x_3): \mathcal{D} \to \R^3$ \'e um gr\'afico e, conseq\"{u}entemente, \'e um mergulho.
\\

Assim, as reflex\~oes no dom\ih nio fundamental sob a superf\ih cie geram a pe\c ca fundamental $P,$ que n\~ao possui auto-interse\c c\~oes e, as interse\c c\~oes de $M$ s\'o ocorrem nos fins catenoidais.
\eject
\ \\
{\huge\bf Refer\^encias}
\ \\
\\
$[1]$ {\it H. Karcher}, Construction of minimal surfaces, Surveys in Geometry, University of Tokyo 1--96, 1989 and Lecture Notes {\bf 12}, SFB256, Bonn, 1989.
\\
$[2]$ {\it F.J. L\'opez} and {\it F. Mart\'\i n}, Complete minimal surfaces in $\R^3$, Publicacions Matematiques {\bf 43} (1999), 341--449.
\\
$[3]$ {\it J.C.C. Nitsche}, Lectures on minimal surfaces, Cambridge University Press, Cambridge, 1989.
\\
$[4]$ {\it R. Osserman}, A survey of minimal surfaces, Dover, New York, 2nd ed, 1986.
\end{document}